\documentclass[12pt]{amsart}

\usepackage{color}
\usepackage{amsbsy}
\usepackage{amssymb,amsmath,amstext,mathrsfs}

\makeatletter
\def\verbatim{\interlinepenalty\@M \@verbatim
  \leftskip\@totalleftmargin\advance\leftskip2pc
  \frenchspacing\@vobeyspaces \@xverbatim}
\makeatother
\newtheorem{thm}{Theorem}[section]
\newtheorem{cor}[thm]{Corollary}
\newtheorem{lem}[thm]{Lemma}
\newtheorem{prop}[thm]{Proposition}

\theoremstyle{definition}
\newtheorem{defn}{Definition}[section]

\theoremstyle{remark}
\newtheorem{rem}{Remark}[section]

\numberwithin{equation}{section}

\newtheorem{ex}[thm]{Example}

\newcommand{\hu}{\hspace{.02in}}
\newcommand{\ca}{\mathrm{card} \hu}

\newcommand{\va}{\left|}
\newcommand{\vb}{\right|}
\newcommand{\vd}{\right\|}
\newcommand{\vc}{\left\|}
\newcommand{\pl}{\left(}
\newcommand{\pr}{\right)}
\newcommand{\tl}{\left\{}
\newcommand{\tr}{\right\}}
\newcommand{\ql}{\left[}
\newcommand{\qr}{\right]}

\newcommand{\ra}{\rightarrow}

\begin{document}
\title[type $\mathrm{I}$ isometric  shifts]{Examples  and counterexamples of type $\mathrm{I}$ isometric shifts}

\author{Jes\'us Araujo}
\address[]{Departamento de Matem\'aticas, Estad\'{\i}stica y Computaci\'on\\Facultad de Ciencias\\
Universidad de Cantabria\\  Avda.
de los Castros, s. n.\\ E-39071 Santander, Spain}
\email[]{araujoj@unican.es}

\thanks{Research partially
supported by the
Spanish Ministry of Science and Education (Grant number MTM2006-14786).}


\maketitle

\begin{abstract}
We provide examples of nonseparable spaces $X$ for which $C(X)$ admits an isometric shift, which  solves in the negative a problem proposed by Gutek {\em et al.} (J. Funct. Anal. {\bf 101} (1991), 97-119).  We also give two independent methods for obtaining separable examples. The first one allows us in particular to construct examples with infinitely many nonhomeomorphic components in a subset of the Hilbert space $\ell^2$. The second one  applies for instance  to sequences adjoined to  any $n$-dimensional compact manifold (for $n \ge 2$) or to the Sierpi\'nski curve. The combination of both techniques leads to different examples involving a convergent sequence adjoined to the Cantor set: one method for the case when the sequence converges to a point in the Cantor set, and the other one for the case when it converges outside.
\end{abstract}

\section{Introduction}
The usual concept of  shift operator in the Hilbert space $\ell^2$ has been introduced in the more general context
of Banach spaces in the following way (see \cite{C, Hl}):   Given a Banach space $E$ over $\mathbb{K}$ (the field of real or complex numbers), a linear operator $T: E \rightarrow E$ is said to be an {\em isometric
shift} if
\begin{enumerate}
\item $T$ is an isometry,
\item The codimension of $T(E)$ in $E$ is $1$,
\item $\bigcap ^{\infty}_{n=1} T^{n}(E)=\{\mathbf{0}\}.$
\end{enumerate}

One of the main settings where isometric shifts have been studied is  $E= C(X)$, that is, the Banach space of 
all $\Bbb{K}$-valued continuous functions defined on a compact and Hausdorff space $X$,  equipped with its usual supremum norm. In this setting, major breakthroughs were made in \cite{GHJR} and \cite{Ha}. On the one hand, in  \cite{GHJR}, Gutek, Hart, Jamison, and Rajagopalan
studied in depth these operators. In particular, using  the well-known Holszty\'{n}ski's Theorem (\cite{Hz}), they  classified them into two types, called type $\mathrm{I}$ and type $\mathrm{II}$. On the other hand,  in \cite{Ha}, Haydon  showed a general method for providing isometric
shifts of type $\mathrm{II}$, as well as concrete examples. 

However, the picture is far from being complete, mainly for two reasons. First, a very basic question has remained open since the publication in 1991 of the seminal paper \cite{GHJR}:
If  $C(X)$ admits an isometric shift, must $X$ be separable?
This question is only meaningful for type $\mathrm{I}$ isometric shifts since it was already proved in \cite[Corollary 2.2]{GHJR} that type $\mathrm{II}$
isometric shifts yield the separability of $X$. Second, it is remarkable the scarcity of examples of isometric
shifts of type $\mathrm{I}$.

Let us recall the definitions. If $T: C(X) \ra C(X)$ is an isometric shift, then there exist a closed subset $Y \subset X$, a continuous and surjective map $\phi: Y \ra X$, and a
function $a\in C(Y)$, $\left|a \right| \equiv 1$, such that
$(Tf)(x)=a(x)\cdot f(\phi(x))$
for all $x\in Y$ and all $f \in C(X)$.
$T$ is said to be of type $\mathrm{I}$ if $Y$ can be taken to be equal to $X \setminus \{p\}$, where $p \in X$ is an isolated point, and is said to be of type $\mathrm{II}$ if $Y$ can be taken equal to $X$. Moreover, if $T$ is of type $\mathrm{I}$, then the map $\phi : X \setminus \{p\} \ra X$ is indeed a  homeomorphism.

All the examples of type $\mathrm{I}$ isometric
shifts provided in \cite{GHJR} are {\em  primitive}, 
that is, they satisfy the following property:
the set $\mathcal{N}:=\{\mathbf{1},\mathbf{2} , \mathbf{3}, ...\}$ is dense in $X$ (where by $\mathbf{n}$ we denote the isolated point $\phi^{n-1} (p)$, for every $n \in \mathbb{N}$). Farid and
Varadarajan (\cite{FV}) provided the first example of a
non-primitive isometric shift, namely $X\setminus {\rm
cl}_{X}(\mathcal{N})$ is a finite nonempty subset (where ${\rm cl}_{X}(\mathcal{N})$ stands for the closure of $\mathcal{N}$ in $X$). In \cite{GN}, the authors
gave an example for which $X\setminus {\rm cl}_{X}(\mathcal{N})$ is the compactification of a
countably infinite set of isolated points and wondered, related to
the non-separability question mentioned above, how big can
$X\setminus {\rm cl}_{X}(\mathcal{N})$ be.

Nevertheless, all these examples of compact spaces which admit type $\mathrm{I}$
isometric shifts have a dense set of isolated points; that is,
they are compactifications of integers. Two very different examples where the set of isolated points
is {\em not} dense were given in \cite{RS4} and \cite{AF3}: The example  in \cite{RS4} involves the Cantor set and is totally disconnected, whereas that in \cite{AF3}  is not. Consequently, the latter example  also gives a negative answer  to
the question addressed in \cite{Hl}, where it was conjectured that the space $X$ cannot have an infinite 
connected component (see also \cite[Corollary 2.1]{GHJR} for an  example involving a primitive isometric shift).
Related to this, one of the main results in \cite{GHJR} states that $C(X)$ does not admit any isometric shifts, whenever $X$ has a countably infinite number of components, all of whom are infinite.

 In this paper, we  will see that the setting of type $\mathrm{I}$ isometric shifts is indeed much more complicated and rich than one could think in principle, and certainly very different from  that of type $\mathrm{II}$ shifts.

In particular, we will give an answer in the negative to the separability question: There are indeed examples of isometric shifts on $C(X)$, with $X$ not separable, and even having $2^{\mathfrak{c}}$ infinite components (see  Sections~\ref{comidajo} and \ref{lucia10-m-07-rememsoper}). On the other hand, we will give two different methods for constructing  isometric shifts of type $\mathrm{I}$ on separable spaces which are {\em not} compactifications of integers. The first one allows us, in some cases, 
to provide examples in which $X \setminus {\rm cl}_{X}(\mathcal{N})$ has a countably infinite number of components, all of whom are infinite,  even in $\mathbb{C}^n$ or in $\ell^2$ (see Theorems~\ref{susovict}, \ref{13dic22}, \ref{nesko}, Remark~\ref{ruido},  and Corollary~\ref{kolya}). The second one deals with the possibility of taking the weight function $a$ constantly equal to $1$, that is, the possibility of finding an isometric shift which is also an algebra isomorphism when considering the restrictions of the images to $C (X \setminus \{\mathbf{1}\})$. The approach must
be necessarily different, and no result related to the first method can be  given (see Proposition~\ref{nde} and 
Theorem~\ref{graciass0}). As particular spaces where this technique can be applied, we have that if $X$ is a convergent sequence 
adjoined  to a compact $n$-manifold,  $n \ge 2$, then $C(X)$ always admits many isometric shifts of this kind (see Example~\ref{delon}). Compare this with the fact that if $M$ is a compact manifold, then $C(M)$ does not admit any isometric shift, as proved in \cite[Theorem 6.1]{FV}. Finally, a combination of both methods allows us to provide different examples involving a sequence adjoined to the Cantor set (see Theorem~\ref{toallas}).

Let us say that the above classification given in \cite{GHJR} is not mutually exclusive, that is, there are isometric 
shifts which are both of type $\mathrm{I}$ and type $\mathrm{II}$. From now on, we will make a distinction between those which are {\em purely} of any of these two types. 
We will say that an isometric shift is of 
type $\mathrm{I}_0$  if it is  {\em not} of type $\mathrm{II}$, and that it is 
 of type $\mathrm{II}_0$ if it is {\em not} of type $\mathrm{I}$.

Some other papers have recently studied questions related to isometric shifts (also defined on other  spaces of functions). Among them, we will mention for instance \cite{AF2, CJW, F, I, JW, RRS, RS1, RS3}, and \cite{RS5} (see also references therein).

\section{Preliminaries and notation}

Throughout topological spaces are assumed to be Tychonoff (which implies in particular they are Hausdorff), but not necessarily compact unless stated otherwise.
The reason for this is that at some points we may ask for weaker conditions so as to be able to produce different compactifications to which results apply. 

Two special spaces we will use are the Cantor set  and the unit circle in $\mathbb{C}$, which will be denoted by $\mathbf{K}$ and $\mathbb{T}$, respectively. Some powers of $\mathbb{T}$ will appear; in particular, $\mathbb{T}^0$ will be  the set  $\{0\}$.

$L^{\infty} (\mathbb{T})$ will be the space of all Lebesgue-measurable essentially bounded {\em complex}-valued functions
on $\mathbb{T}$, and  $\mathfrak{M}$ will be its maximal ideal space. $m$ will denote the Lebesgue measure on $\mathbb{T}$.

In general,  given a continuous map $f$ defined on a space $X$, we also denote by $f$ its restrictions to subspaces of $X$ and its extensions to other spaces  containing $X$. Given a (surjective) homeomorphism $\phi : X \ra X$, we denote $\phi^n := \underbrace{\phi  \circ \cdots \circ \phi}_n$, for $n \in \mathbb{N}$.

For a set $A$, $\ca A$ denotes its  cardinal if $A$ is finite, and $\ca A := \infty$ otherwise. If there is not possibility of confusion, given a set of indexes $A$ and spaces $X_s$ all of them equal to $X$, we denote the product  $\prod_{s \in A} X_s$ as $X^A$ and $X^{\kappa}$, indistinctly, where $\kappa$ is the cardinal of $A$. As usual,  $\aleph_0$ and $\mathfrak{c}$ denote the cardinals of $\mathbb{N}$ and $\mathbb{R}$, respectively.
A result which will be widely used is the following:  Given a separable space $W$ with at least two points and a cardinal $\kappa$, the power $W^{\kappa}$ is separable if and only if $\kappa \le \mathfrak{c}$ (see \cite[Theorem 16.4]{W}). We do not need to make
any special assumptions about cardinals between $\aleph_0$ and $\mathfrak{c}$, and results will be valid if we assume whether or  not the 
Continuum Hypothesis holds.

We will often refer to the topological sum (or disjoint union) of a family $\pl X_{\alpha} \pr_{\alpha \in \Lambda}$ of pairwise disjoint topological spaces. Recall that
it is the union $\bigcup_{\alpha \in \Lambda} X_{\alpha}$ endowed with the topology consisting of unions of open subsets of these spaces,
and that in particular each $X_{\alpha}$ becomes  closed and open (see \cite[p. 65]{W}).
If $\Lambda:= \{\alpha_1 ,\ldots, \alpha_n, \ldots\}$ is at most countable, then the topological sum will be denoted  by $\sum_{n} X_{\alpha_{n}}$ or $X_{\alpha_1} + \cdots + X_{\alpha_n} + \cdots $.

For a topological space $X$, we will denote by $\beta X$ its Stone-\v{C}ech compactification. Also,  $\mathcal{N} \cup \{\pmb{\infty}\}$ will be  the one-point compactification of $\mathcal{N}$.

On the other hand, some examples will be given contained in the usual Hilbert space $\ell^2$. We will not use the Hilbert space structure, but its topology. Consequently, we may view it both as a real or complex vector space, since
both are homeomorphic.

Throughout "homeomorphism" will be synonymous with "{\em surjective} homeomorphism".

We will usually write $T = T [a, \phi, \Delta]$ to describe a codimension $1$ linear isometry $T: C(X) \ra C(X)$, where 
$X$ is compact and contains $\mathcal{N}$. It means that $\phi : X \setminus \{\mathbf{1}\} \ra X$ is a homeomorphism, satisfying in particular 
$\phi (\mathbf{n+1}) = \mathbf{n}$ for all $n \in \mathbb{N}$. It also means that $a \in C(X \setminus \{\mathbf{1}\})$, $\va a \vb \equiv 1$, and that $\Delta$ is a continuous linear functional on $C(X)$ with $\vc \Delta \vd \le 1$.
Finally, the description of $T$ we have is  $(Tf) (x) = a(x) f(\phi(x))$,  when $x \neq \mathbf{1}$, and $(Tf) (\mathbf{1}) = \Delta (f)$, for every $f \in C(X)$.

A special element in $C(X)'$ is the evaluation map at a point $x \in X$, which will be denoted by $\Gamma_x$.

All our results will be valid in the real and complex settings, unless otherwise stated. The only exceptions are the following: Results exclusively given for $\mathbb{K} = \mathbb{C}$ appear just in Section~\ref{emphado}.
The only result valid just for the case $\mathbb{K} = \mathbb{R}$ is given in Example~\ref{torero}. $C_{\mathbb{C}} (X)$ and $C_{\mathbb{R}} (X)$    will denote the Banach spaces of functions on $X$ taking complex and real values, respectively.

\medskip

We now give some basic definitions that will be used in particular in Sections~\ref{unya} and \ref{majorca}.

\begin{defn}
A subset $\mathbb{P} \subset \mathbb{N}$ is an {\em initial subset} if $1 \in \mathbb{P}$, and $n-1 \in \mathbb{P}$ whenever $n \ge 2$ and $n \in \mathbb{P}$.
\end{defn}

\begin{defn} 
Given an initial subset $\mathbb{P}$, a   sequence $(p_n)_{n \in \mathbb{P}}$ of natural numbers is said to be $\mathbb{P}$-compatible if 
\begin{itemize}
\item if $\mathbb{P} = \{1\}$, then  $p_1 > 1$, and
\item if $\mathbb{P} \neq \{1\}$, then  $p_{n+1}$ is an even multiple of $p_n$ for every $n$. 
\end{itemize}
\end{defn}

Our first method will allow us to define isometric shifts with special features. We can create spaces subject to certain restrictions concerning denseness, which cannot be found in examples of isometric shifts  of type $\mathrm{II}_0$. For this reason we introduce the following definition. 

\begin{defn}\label{salto}
Let $X$ be compact, and suppose that $T =T[a, \phi, \Delta]:   C(X) \ra C(X)$
is a {\em non-primitive}  isometric shift of type $\mathrm{I}$.  For $n \in \mathbb{N}$, we say that $T$ is {\em
$n$-generated} if $n$ is the least number with the following  property: There exist $n$ points $x_1 , \ldots ,
x_{n} \in X\setminus {\rm cl}_{X}(\mathcal{N})$ such that the set
\[\{\phi^k (x_i) : k \in {{\mathbb{Z}}}, i \in \{1,  \ldots, n\}\}\] is dense in $X\setminus {\rm cl}_{X}(\mathcal{N})$.

We say that $T$ is $\infty$-generated if it is not
$n$-generated for any $n\in \mathbb{N}$.
\end{defn}

Notice that the above definition does not make sense when $T$ is an isometric shift  of type $\mathrm{II}_0$. In fact,  it is proved in  
\cite[Theorem 2.5]{GHJR} that for such isometric shifts,  there exists a point $x \in X$ such that $\{\phi^k (x) : n \in \mathbb{Z}\}$
is dense in $X$. In particular, if $T$ is a non-primitive isometric shift of both types $\mathrm{I}$ and $\mathrm{II}$, then it is $1$-generated; an example was mentioned above, as appearing in \cite{RS4}, and another example can be found in the proof of Theorem~\ref{notatillo} (it is denoted $T_1$ and defined in the complex setting, but is valid in the real case as well).

\section{Nonseparable examples}\label{comidajo}

We start this section with an answer to the problem of the separability question.

\begin{thm}\label{podmo-skilo-11m}
 $C \pl \mathfrak{M} + \mathcal{N} \cup \{\pmb{\infty}\} \pr $ admits an isometric shift.
\end{thm}

Once we have a first example, we can get more. For instance, the next result is essentially different in that it provides examples with $2^{\mathfrak{c}}$ infinite connected components.

\begin{thm}\label{sinregalocumple}
Let $\kappa$ be any cardinal such that $1 \le \kappa \le \mathfrak{c}$. Then  $C \pl \mathfrak{M} \times \mathbb{T}^{\kappa} + \mathcal{N} \cup \{\pmb{\infty}\} \pr$ admits an  isometric shift.
\end{thm}

Finally, we can also give examples with just one infinite component.

\begin{thm}\label{notatillo2}
Let $\kappa$ be any cardinal such that $1 \le \kappa \le \mathfrak{c}$. Then $C \pl \mathfrak{M} + \mathbb{T}^{\kappa} + \mathcal{N} \cup \{\pmb{\infty}\} \pr $ admits  an isometric shift.
\end{thm}

\begin{rem}\label{perrol}
Even if our space $\mathfrak{M}$ is based on an algebra of complex-valued functions, Theorems~\ref{podmo-skilo-11m}, \ref{sinregalocumple}, and \ref{notatillo2} are valid both if $\mathbb{K}= \mathbb{R}$ or $\mathbb{C}$. Nevertheless, there are examples that can be constructed just in the complex setting (see Theorem~\ref{jarboresquil} and Example~\ref{torero}).
\end{rem}

\section{Separable examples}\label{unya}

In this section, we provide some  results which easily lead to different examples. We do it in particular
contexts where our methods can be applied. These contexts involve separately: 1)  copies of (finite or infinite) 
separable powers of $\mathbb{T}$; 2)  copies of separable {\em infinite} powers of {\em any} compact spaces; 3) any compact manifold and some sets as the Sierpi\'nski curve; 4) two different situations (the only possible involving a convergent sequence) for the Cantor set. 

Theorems~\ref{13dic22} and \ref{nesko} are consequences of results given in Section~\ref{majorca} (see also Remarks~\ref{galletuca} to \ref{masiosare}). In both cases it is possible not only to ensure the existence of  one isometric shift, 
but also to give infinitely many of them attending to  the notion of $n$-generated shift.

\begin{ex}\label{gruny}
Let $\mathbb{P} := \{1, \ldots, k\}$, and let $(p_n)_{n \in \mathbb{P}}$ be a $\mathbb{P}$-compatible sequence. 
Let $\{m_1, \ldots, m_k\} \subset \mathbb{N} \cup \{0\}$.
Then  there exists  a $k$-generated isometric shift   of type $\mathrm{I}_0$ on $C(X)$, where $$X = \underbrace{\mathbb{T}^{m_1} + \cdots + \mathbb{T}^{m_1}}_{p_1} + \cdots + \underbrace{\mathbb{T}^{m_k} + \cdots + \mathbb{T}^{m_k}}_{p_k} + \mathcal{N} \cup \{\pmb{\infty}\}.$$ 
\end{ex}
The above example can be generalized in several ways, attending to the powers of $\mathbb{T}$ involved or to the number of copies of each of them.
\begin{ex}\label{pasaprau}
Let $\mathbb{P}$ and $(p_n)_{n \in \mathbb{P}}$ be as in Example~\ref{gruny}. Consider also  $\{m_1, \ldots, m_l\} \subset \mathbb{N} \cup \{0\}$, where 
$l \le k$. Make a partition of $\{p_n : n \in \mathbb{P}\}$  into $l$ subsets $\mathscr{A}_1 , \ldots, \mathscr{A}_l$, and let $q_i$ the sum of the numbers contained in each $\mathscr{A}_i$. 
Then  there exists  a $k$-generated isometric shift   of type $\mathrm{I}_0$ on $C(X)$, where
$$X = \underbrace{\mathbb{T}^{m_1} + \cdots + \mathbb{T}^{m_1}}_{q_1} + \cdots + \underbrace{\mathbb{T}^{m_l} + \cdots + \mathbb{T}^{m_l}}_{q_l} + \mathcal{N} \cup \{\pmb{\infty}\}.$$
\end{ex}

\begin{ex}\label{cafard}
It is also possible to give examples where the number of components is infinite. In this case, we get
examples of isometric shifts which are not finitely generated. These examples are obtained in $\ell^2$, endowed with the norm topology, and obviously
in the construction a process of compactification must be made. More concretely,
let  $(p_n)_{n \in \mathbb{N}}$ be an $\mathbb{N}$-compatible sequence, and let $(q_n)_{n \in \mathbb{N}}$ be as in Example~\ref{pasaprau}, that is, we make a partition $(\mathscr{A}_i)_{i \in \mathbb{N}}$ of $\{p_n: n \in \mathbb{N}\}$ into finite subsets, and define $q_n := \sum_{k \in \mathscr{A}_n} k$ for each $n$. 
Suppose also that $(m_n)_{n \in \mathbb{N}}$ is  a sequence of (pairwise different) elements in $\mathbb{N} \cup \{0\}$.  Then  there exists  an $\infty$-generated isometric shift    on $C(X)$ (obviously  of type $\mathrm{I}_0$), 
where 
 $X \subset \ell^2$  consists of a compactification of 
$$\pl \sum_{n=1}^{\infty} \underbrace{\mathbb{T}^{m_n} + \cdots + \mathbb{T}^{m_n}}_{q_n} \pr + \mathcal{N} \cup \{\pmb{\infty}\}.$$
On the other hand, exactly the same conclusion
can also be reached if one of the $m_n$ is taken to be equal to $ \aleph_0$.
\end{ex}

\begin{ex}\label{luzinfinita}
In the above examples, we are considering sums (or compactifications of sums) of finite copies of different powers of $\mathbb{T}$ (being these powers finite or countable). It is also possible to give similar examples where infinite copies are 
envolved, even without requiring the $m_n$ to be finite, in the following way: Let $(\kappa_n)$ be a {\em finite} or {\em infinite} sequence of (pairwise different) cardinals, each satisfying $0 \le \kappa_n \le \mathfrak{c}$. 
Then  there exists  an $\infty$-generated isometric shift   on $C(X)$, where 
 $X$  consists of a compactification of
 $$\pl \sum_{n} \underbrace{\mathbb{T}^{\kappa_n} + \cdots + \mathbb{T}^{\kappa_n} + \cdots}_{\aleph_0} \pr + \mathcal{N} \cup \{\pmb{\infty}\}.$$
\end{ex}

The combination of the ideas of obtaining finite copies of some powers and infinite copies of some other, as
well as the previous examples, can be summarized in the following result.

\begin{thm}\label{13dic22}
Let $\mathbb{M}$ and $ \mathbb{P}$ be initial subsets of $\mathbb{N}$ with $\mathbb{M} \subset \mathbb{P}$. Let  
$(\kappa_n)_{n \in \mathbb{M}}$ be  a sequence of pairwise different  cardinals with $0 \le \kappa_n \le \mathfrak{c}$, and let 
$(p_n)_{n \in \mathbb{P}}$ be a $\mathbb{P}$-compatible sequence. Suppose that $(\mathscr{A}_i)_{i \in \mathbb{M}}$ is a partition of $\{p_n : n \in \mathbb{P}\}$, and that 
$q_n := \sum_{k \in \mathscr{A}_n} k$ for every $n \in \mathbb{M}$.

Then, for $N= \ca \mathbb{P}$, there exist a complex dual Banach space $\mathbf{B}'$ (endowed with the $\mathrm{weak}^*$-topology) and a compactification 
 $X \subset \mathbf{B}'$  of   
 $$\pl \sum_{n\in \mathbb{M}} \underbrace{\mathbb{T}^{\kappa_n} + \cdots + \mathbb{T}^{\kappa_n}}_{q_n} \pr + \mathcal{N} \cup \{\pmb{\infty}\}$$
such that $C(X)$ admits  an $N$-generated isometric shift   of type $\mathrm{I}_0$.
Moreover, if $\mathrm{s}:= \sup \{\kappa_n : n \in \mathbb{M}\} < \aleph_0$, then $\mathbf{B}'$ may be taken to be  $\mathbb{C}^{\mathrm{s}}$. Similarly, if $\mathrm{s} = \aleph_0$, then $\mathbf{B}'$ may be taken to be   $\ell^2$ with the norm topology.
\end{thm}

A similar result can be given involving {\em infinite} powers of separable spaces, as for example the Cantor set 
(homeomorphic to $\{0,1\}^{\mathbb{Z}}$), $\{0,1\}^{\mathbb{R}}$,
 $[0,1]^{\mathbb{Z}}$ or $[0,1]^{\mathbb{R}}$, among others.

\begin{thm}\label{nesko}
Let $\mathbb{M}$, $\mathbb{P}$, $(p_n)_{n \in \mathbb{P}}$,  and $(q_n)_{n \in \mathbb{M}}$ be as in Theorem~\ref{13dic22}.  For each $n \in \mathbb{M}$, let $K_n$ and $\kappa_n$ be a  separable (nonempty) compact space and a cardinal with $\aleph_0 \le \kappa_n \le \mathfrak{c}$, respectively.
  Assume also that $K_n^{\kappa_n} \neq K_k^{\kappa_k}$ whenever $n \neq k$, $n, k \in \mathbb{M}$.

Then, for $N= \ca \mathbb{P}$, there exists  a compactification 
 $X$  of 
 $$\pl \sum_{n\in \mathbb{M}} \underbrace{K_n^{\kappa_n} + \cdots + K_n^{\kappa_n}}_{q_n} \pr + \mathcal{N} \cup \{\pmb{\infty}\}$$ such that $C(X)$ admits  an $N$-generated isometric shift   of type $\mathrm{I}_0$. 
Moreover, if $K_n$ is metrizable and $\kappa_n = \aleph_0$   for every $n \in \mathbb{M}$, then $X$ may be taken to be contained in  $\ell^2$, endowed  with the norm topology.
\end{thm}

\begin{rem}\label{cinturone}
Recall that in Example~\ref{cafard}, we are assuming that all sets $\mathscr{A}_i$ are  finite. In principle, in Theorems~\ref{13dic22} and \ref{nesko}, each set $\mathscr{A}_i$ can be either finite or infinite. Of course, if some $\mathscr{A}_i$ is infinite, then  $q_i = \infty$. 
\end{rem}

\begin{rem}\label{galletuca}
In Theorem~\ref{nesko}, we are not assuming that spaces $K_n$ are necessarily pairwise different. In the same way, in Theorems~\ref{13dic22} and \ref{nesko}, cardinals $\kappa_n$ are assumed to be infinite but not necessarily different. Of course there are at least two possible cardinals like this, namely  $\aleph_0$ and  $\mathfrak{c}$. On the one  hand, if we assume that 
 the Continuum Hypothesis holds, then these are the only cardinals satisfying our requirements. In that case,  $\mathbb{M}$ is not necessarily finite in Theorem~\ref{nesko}, since we could be dealing with powers of very different spaces
 $K_n$, so it could be the case that $\mathbb{M} = \mathbb{N}$ and $\kappa_n$ is the same cardinal for every $n$. 
 
 On the other hand, if we assume that Continuum Hypothesis does not hold, then we may have a variety of different cardinals satisfying our assumptions, and each one can provide a different power of $\mathbb{T}$ in Theorem~\ref{13dic22}.
 In particular, if we assume that all $\kappa_n$ are pairwise different, then all spaces $K_n$ may be taken to be the same in Theorem~\ref{nesko}, as long as $\ca K_n >1$.
\end{rem}

\begin{rem}\label{dia6-4-07}
In Theorems~\ref{13dic22} and \ref{nesko}, metrizability of the examples is achieved under some conditions. Indeed metrizability is sometimes a consequence of the following fact: If $C(X)$ admits an isometric shift, then $X$ must satisfy the countable chain condition (see \cite[Theorem 1.4]{GN},  \cite[Theorem 7]{RS4} or  \cite[Lemma 5.6]{RS3}). For this reason, if an example $X$ is  a weakly compact subset of a Banach space, then it is also metrizable (see \cite[Proof of Corollary 4.6]{Ro}, and see also \cite{ARN} and \cite{M} for closely related results). This implies that we cannot in general replace in  Theorem~\ref{13dic22} the statement referring to the $\mathrm{weak}^*$-topology with another one referring to the $\mathrm{weak}$-topology.
\end{rem}

\begin{rem}\label{inca}
In Theorem~\ref{nesko}, it is possible to include $\mathbb{T}$ as one (or several) of the spaces $K_n$, but the construction of the associated homeomorphism will be different from the one we use in Theorem~\ref{13dic22}.
\end{rem}

\begin{rem}\label{ruido}
Notice that each copy of a power of $\mathbb{T}$ given in Theorem~\ref{13dic22} is a connected component of 
$X\setminus {\rm cl}_{X}(\mathcal{N})$, and they are indeed all the infinite connected components of $X$ (if we assume  $\kappa_n \neq 0$ for every $n$). The same
comment applies to copies of powers of $K_n$ in Theorem~\ref{nesko} when all $K_n$ are connected. This implies that we can
construct examples where we may decide at will on the number of (different) infinite connected components of $X\setminus {\rm cl}_{X}(\mathcal{N})$.
\end{rem}

\begin{rem}\label{masiosare}
In Theorems~\ref{13dic22} and \ref{nesko}, we deal with powers of $\mathbb{T}$ and spaces $K_n$, which are compact, and we refer to "a compactification". Similar results can be given to powers of noncompact spaces, and then work with a compactification (as will be stated in Section~\ref{majorca}). It will be always possible to do it with the Stone-\v{C}ech compactification, but many times it will even be possible to work with an unrelated one, as we will see.
\end{rem}

All the previous  results are  based on our first method. As for the second method, given in Section~\ref{luzia}, it provides some other examples.
\begin{ex}\label{delon}
In some cases, given a compact space $W$, it is possible to find a compact space $X$ such that $C(X)$ admits an isometric shift  of type $\mathrm{I}_0$ for which $a \equiv 1$, and such that $W = X \setminus \mathcal{N}$. This happens for instance in any of the three following cases:

\begin{itemize}
\item if $W$ is a separable infinite power of a compact space with at least two points,
\item if $W$ is a compact  $n$-manifold (with or without boundary), for $n \ge 2$,
\item if $W$ is the Sierpi\'nski curve.
\end{itemize}
\end{ex}

Also in the case when the weight $a$ is equal to $1$, it is possible to find examples where the number of infinite connected components is $n$, for every $n \in \mathbb{N}$ (see Corollary~\ref{labaderu}).

\medskip

We finally give a general theorem involving the Cantor set. The set $\{x_n : n \in \mathbb{N}\}$ will play the r\^{o}le of $\mathcal{N}$. To get it we combine the results of Sections~\ref{majorca} and \ref{luzia}.
\begin{thm}\label{toallas}
Suppose that  $(x_n)$ is a nonconstant sequence in $\mathbb{R} \setminus \mathbf{K}$ which converges to a point $L \in \mathbb{R}$. Let $X := \mathbf{K} \cup \{x_n : n \in \mathbb{N}\} \cup \{L\}$. We have
\begin{itemize}
\item If $L \notin \mathbf{K}$, then for each $n \in \mathbb{N} \cup \{\infty\}$, there exists an isometric shift  of type $\mathrm{I}_0$ on $C(X)$ which is $n$-generated.
\item If $L \in \mathbf{K}$, then there exists an isometric shift  of type $\mathrm{I}_0$ on $C(X)$ which also satisfies the additional condition that $a \equiv 1$.
\end{itemize}
\end{thm}

\begin{rem}
Recall that, as mentioned by the authors, by \cite[Theorem 1.9]{GN} we can conclude that if $X$ consists of a convergent sequence adjoined to  a non-separable Cantor cube, then $C(X)$ does not admit an isometric shift. This is not true in the separable case,  as shown in \cite[Example 20]{RS4} for an isometric shift of both types $\mathrm{I}$ and $\mathrm{II}$. Example~\ref{delon} and Theorem~\ref{toallas} say also the contrary  in the separable case for isometric shifts which are not of type $\mathrm{II}$.  Theorem~\ref{toallas} provides indeed  completely different families of isometric shifts  of type $\mathrm{I}_0$.   Recall finally that in \cite[Theorem 1]{Ha}, it is proved that in the case where $\mathbb{K} = \mathbb{R}$, $C(\mathbf{K})$
admits an isometric shift (necessarily  of type $\mathrm{II}_0$ because $\mathbf{K}$ has no isolated points).
\end{rem}

\section{Some results on product spaces and $\mathbb{L}$-transitivity}

Recall that given a  topological space $W$ and a homeomorphism $\phi$  from $W$ onto itself, the semiflow $(W, \phi)$
is defined as the sequence $( \phi^n)_{n \in \mathbb{N}}$ of iterates of $\phi$. Sometimes we will say that $\phi$ is itself a semiflow, if it is clear to which topological space we are referring to. On the other hand, for a  point $x \in W$, we denote $\mathrm{Orb}^+ (\phi, x) := \{\phi^n (x) :n \in \mathbb{N}\}$. Our definitions will refer to this semiorbit.

\begin{defn}
We say that a semiflow $(W, \phi)$ is 
\begin{itemize}
\item (topologically) {\em transitive} if there exists a point $w \in W$  
such that $\mathrm{Orb}^+ (\phi, w)$  is dense in $W$.
\item (topologically) {\em totally transitive} if there exists a point $w \in W$  
such that $\mathrm{Orb}^+ (\phi^k, w)$  is dense in $W$ for every $k \in \mathbb{N}$.
\end{itemize}
\end{defn}
We will omit the adverb {\em topologically} when referring to these concepts. In particular, in the above contexts,
we will also say that $(W, \phi, w)$ is transitive and totally transitive, respectively. If there is no possibility of confusion we may also say that $W$ or that $\phi$ is (totally) transitive.
 
\begin{rem}\label{trankizia}
In this paper, we will use the {\em above} definitions of transitivity and total transitivity. Nevertheless  we must 
remark that there are some others that we will not use: In particular, $\phi$ is sometimes called transitive if there 
exists a point $x \in W$ such that $\{\phi^n (x) :n \in \mathbb{Z}\}$ is dense in $W$, and totally transitive if 
$\phi^n$ is transitive (in the latter sense) for every $n \in \mathbb{N}$. Both definitions of transitivity agree when 
$W$ is a complete separable metric space without isolated points: Indeed, in such case the set of points satisfying that $\mathrm{Orb}^+ (\phi, w)$  is 
dense in $W$ is a dense $G_{\delta}$-set (see \cite[Chapter 18]{O}). This implies also that the set of points $w$ satisfying that $\mathrm{Orb}^+ (\phi^n , w)$ is dense for every $n$ is also a dense $G_{\delta}$-set,  
and both definitions of total transitivity are equivalent.
\end{rem}

All our results on isometric shifts will depend on some of the above notions: Plain transitivity will be necessary
to get results in Section~\ref{luzia}. Most separable and nonseparable examples will rely on a flow for which a given power is transitive. This will allow in particular to obtain $n$-generated isometric shifts, for $n \neq \infty$. On the other hand, we will use total transitivity to construct $\infty$-generated isometric shifts in the {\em separable} case. To do so, we still need  another definition, that of $\mathbb{L}$-transitivity, which will be introduced later.

\medskip

Before giving an example  which will be used to provide examples of isometric shifts,
let us recall that
given a family of semiflows $\pl \pl Z_i , h_i \pr \pr_{i \in I}$, it is possible to define a new semiflow $(h_i)_{i \in I}: \prod_{i \in I} Z_i \ra \prod_{i \in I} Z_i$ coordinatewise, that is, $(h_i)_{i \in I} \pl (z_i)_{i \in I}) :=  (h_i (z_i) \pr_{i \in I}$ for each point $(z_i)_{i \in I}$.  The new semiflow $\pl \prod_{i \in I} Z_i, (h_i)_{i \in I} \pr$ 
will be called  product of $\pl Z_i , h_i \pr$, $i \in I$.
\begin{ex}\label{dominicus}
If we take the unit circle $\mathbb{T} $, and $\rho$ is an irrational number, then the rotation 
semiflow $[\rho] : \mathbb{T} \ra \mathbb{T}$ sending each $z \in \mathbb{T}$ to $z e^{  2 \pi \rho i}$ satisfies that 
$(\mathbb{T}, [\rho], z)$ is transitive for every 
$z \in \mathbb{T}$ (see \cite[Proposition III.1.4]{V}).  Indeed, since the same applies to $k \rho$ for every $k \in \mathbb{N}$, we have that $(\mathbb{T}, [\rho], z)$ is totally transitive for every 
$z \in \mathbb{T}$. It is easy to see that this fact can be generalized to separable powers of $\mathbb{T}$ (similarly as it is mentioned for finite powers in \cite[III.1.14]{V}):  
Let $\Lambda := \{\rho_{\alpha} : \alpha \in \mathbb{R}\}$ be a set of irrational numbers linearly independent over
 $\mathbb{Q}$. If $\mathbb{P}$ is any nonempty subset of $\mathbb{R}$ and $\pl \mathbb{T}^{\mathbb{P}}, \ql \rho_{\alpha} \qr_{\alpha \in \mathbb{P}}  \pr$ is the product of rotation  semiflows 
$\pl \mathbb{T}_{\alpha} , [\rho_{\alpha}]\pr$, then   
$\pl \mathbb{T}^{\mathbb{P}}, \ql \rho_{\alpha} \qr_{\alpha \in \mathbb{P}}, \pl z_{\alpha}\pr_{\alpha \in \mathbb{P}} \pr $ is totally transitive for any choice of  points $z_{\alpha} \in \mathbb{T}_{\alpha}$, $\alpha \in \mathbb{P}$.
\end{ex}
  
The following lemma is well known and easy to prove. 

\begin{lem}\label{aqua}
Suppose that $W, Z$ are topological spaces, and that $\phi :  W \ra W$ and $\psi: Z \ra W$ are   homeomorphisms. If $(W, \phi, w)$ is transitive (respectively, totally transitive), then $(Z, \psi^{-1} \circ \phi \circ \psi, \psi^{-1} (w))$ is transitive (respectively, totally transitive). 
\end{lem}

The above lemma implies in particular that we can work with spaces homeomorphic to those we want to get results on.
This will be useful in some instances.

\begin{thm}\label{tabaclaro-v59}
Let  $W$ be a separable topological   space homeomorphic to $W^{\mathbb{N}}$. Then there exists a   homeomorphism $\phi :W \ra W$ which fixes a point of $W$ and has a dense set of periodic points, and  such that $(W, \phi)$ is totally  transitive.
\end{thm}

\begin{proof}
We will prove  that $(W^{\mathbb{Z}}, \Sigma)$ is totally transitive, where $\Sigma$ is the usual {\em bilateral shift map} on $W^{\mathbb{Z} }$ (defined as 
$(\Sigma x)_m = x_{m+1}$ for every $x = \pl x_n \pr_{n \in \mathbb{Z} } \in W^{\mathbb{Z}}$ and $m \in \mathbb{Z}$). 

Suppose  that $\{x_1, x_2, \ldots, x_n, \ldots \}$  is a dense subset of $W$. We assume that this set is infinite, 
but a similar proof would work in the finite case. For each $n \in \mathbb{N}$, consider $D_n := \{x_1,  \ldots , x_n\}$ and   the Cartesian product $D_n^n$, which consists of $n$-tuples $P_i^n$ ($i=1, \ldots, n^n$) of the form   \[P_i^n = (P_i^n (1), P_i^n (2), \ldots,  P_i^n (n)),\] with each $P_i^n (j) \in D_n$.

We are going to define inductively a special element  $w = (w_n) \in W^{\mathbb{Z}}$. 
Let $A:= \{(i, j) \in \mathbb{N}^2 : i \le j^j\}$,
and let $\psi :  \mathbb{N} \ra \mathbb{Z} \times \mathbb{N} \times A$ be  a bijective map. For each $n \in \mathbb{N}$, denote by $\psi_1 (n)$ and $\psi_2 (n)$ the first and  second coordinates of $\psi (n)$, respectively. 
As for the third coordinate of $\psi (n)$, which  is an element of $A$,  we denote it by $(\psi_3 (n), \psi_4 (n))$. We start the induction process 
by considering $\psi_1 (1) \in \mathbb{Z}$ and $\psi_2 (1) \in \mathbb{N}$. Pick $n_1 \in \mathbb{N}$ of the form $n_1 =   \psi_1 (1) + l_1 \psi_2 (1) $, for $l_1 \in \mathbb{N}$, 
and define $w_{n_1 + 1} := P_{\psi_3 (1)}^{\psi_4 (1)} (1)$,  $w_{n_1 + 2} := P_{\psi_3 (1)}^{\psi_4 (1)} (2)$,
 \ldots, $w_{n_1 + \psi_4 (1)} := P_{\psi_3 (1)}^{\psi_4 (1)} (\psi_4 (1))$. Next suppose that we have defined
 $n_1 ,  \ldots, n_k$ (of the form $\psi_1 (1) + l_1 \psi_2 (1), \ldots, \psi_1 (k) + l_k \psi_2 (k)$, respectively) in such a way that 
$n_i + \psi_4 (i) \le n_{i+1}$ for $i=1, \ldots, k-1$. Then we take   $n_{k +1} := \psi_1 (k+1) + l_{k+1} \psi_2 (k+1)$
for $l_{k+1} \in \mathbb{N}$, 
such that $n_k + \psi_4 (k) \le n_{k+1}$, and define
 $w_{n_{k+1} + 1} := P_{\psi_3 (k+1)}^{\psi_4 (k+1)} (1)$, 
 $w_{n_{k+1} + 2} := P_{\psi_3 (k+1)}^{\psi_4 (k+1)} (2)$, \ldots,
 $w_{n_{k+1} + \psi_4 (k+1)} := P_{\psi_3 (k+1)}^{\psi_4 (k+1)} (\psi_4 (k+1))$.

Indeed the rest of coordinates  of $w$, that is, those which cannot be obtained following this process, can be taken without any requirements. For this reason we fix a point $x \in W$, and define $w_n := x$ for every $n \notin \bigcup_{k=1}^{\infty} \{n_k +1 , \ldots, n_k + \psi_4 (k) \}$. 

It is  now easy to check that  $(W^{\mathbb{Z}}, \Sigma, w)$ is totally transitive. 
Also we have that $\Sigma (\mathbf{0}) = \mathbf{0}$ (being $\mathbf{0}$ the point with all coordinates equal to $0$). 
Finally, we get the desired conclusion from Lemma~\ref{aqua}.
\end{proof}

Taking into account the properties of the homeomorphism  given in Theorem~\ref{tabaclaro-v59}, we immediately obtain  the following corollary.

\begin{cor}\label{purusa}
Let  $W$ be a separable topological space, and let $\kappa $ be an infinite cardinal, $\kappa \le \mathfrak{c}$. Then there exists a homeomorphism $\phi$ from $W^{\kappa}$ onto itself, having a fixed point and a dense set of periodic points, such that $\pl W^{\kappa}, \phi \pr$ is totally transitive.
\end{cor}

\begin{rem}
Example~\ref{dominicus} and Corollary~\ref{purusa} give us two different ways for  constructing totally transitive semiflows
on  infinite powers $\mathbb{T}^{\mathbb{N}}$, $\mathbb{T}^{\mathbb{R}}$, and more in general $\mathbb{T}^{\kappa}$,
if $\kappa$ is the cardinal of any infinite subset of $\mathbb{R}$. Both ways will be used to obtain different kinds of examples of isometric shifts  of type $\mathrm{I}_0$.
\end{rem}

Finally, by Lemma~\ref{aqua}, the fact that the Cantor set is homeomorphic to $\{0,1\}^{\mathbb{N}}$ gives us the following corollary.
\begin{cor}\label{canutoru}
There exists a homeomorphism $\phi$ from $\mathbf{K}$ onto itself, having a fixed point and a dense set of periodic points, such that $(\mathbf{K}, \phi)$ is  totally transitive.
\end{cor}

Corollary~\ref{canutoru} will be used to construct examples of isometric shifts  according to the method
given in Section~\ref{majorca}, where the following definition will be fundamental.

\begin{defn}
Let $\mathbb{P}$ and $ \mathbb{L}$ be nonempty subsets of $\mathbb{N}$. For each $n \in \mathbb{P}$, let $(Z_n, h_n)$ be a semiflow, and let $1_n \in Z_n$ be such that $(Z_n, h_n, 1_n)$ is transitive.  We say that 
the (finite or infinite) sequence $((Z_n, h_n, 1_n))_{n \in \mathbb{P}}$ is $\mathbb{L}$-{\em transitive} if for every $k \in \mathbb{L}$ and $i \in \mathbb{N}$, the point $(h_n^i(1_n))_{n \in \mathbb{P}} \in \prod_{n \in \mathbb{P}} Z_n$  belongs to the closure of $\mathrm{Orb}^+ ((h_n^k)_{n \in \mathbb{P}}, (1_n)_{n \in \mathbb{P}})$.
\end{defn}

We will write $(Z_n, h_n, 1_n)_{n \in \mathbb{P}}$ for short, or $(Z_n, h_n, 1_n)$ if there is no confusion about the set $\mathbb{P}$.

\begin{ex}
Let $(Z, h, 1)$ be transitive. Obviously, if $\mathbb{P}$ is a nonempty subset of $\mathbb{N}$ and we take $(Z_n, h_n, 1_n) = (Z, h, 1)$ for every $n \in \mathbb{P}$, then 
$\mathscr{A}:= (Z_n, h_n, 1_n)_{n \in \mathbb{P} }$ is $\mathbb{L}$-transitive, for $\mathbb{L} := \{1\}$. In the same
 way, it is also immediate that if $\mathbb{L} $ is any nonempty subset of $\mathbb{N}$ and $(Z, h, 1)$ is $\mathbb{L}$-transitive, then $\mathscr{A}$ is $\mathbb{L}$-transitive. With this example we also see that our notion of 
$\mathbb{L}$-transitivity does not imply  that $\pl \prod_{n \in \mathbb{P}} Z_n, (h_n), (1_n) \pr$ is transitive.
\end{ex}

\begin{ex}\label{zarrau}
Let $\Lambda$ and 
$(\mathbb{T}_{\alpha} , [\rho_{\alpha}])$ be as in Example~\ref{dominicus}. If  $\mathbb{P} \subset \mathbb{N}$ and  $\{\mathbb{P}_n : n \in \mathbb{P}\}$ is a pairwise disjoint family of (nonempty) subsets of $\mathbb{R}$, then the
sequence $ \pl \mathbb{T}^{\mathbb{P}_n}, \ql \rho_{\alpha} \qr_{\alpha \in \mathbb{P}_n}, \pl z_{\alpha}\pr_{\alpha \in \mathbb{P}_n} \pr_{n \in \mathbb{P}}$ 
is
$\mathbb{N}$-transitive for any choice of points $z_{\alpha}$.

Consequently we see that, if $(n_i)$ and $(l_j)$ are any (finite or infinite) sequences in $\mathbb{N} \cup \{\aleph_0\}$, then it is possible to 
construct an $\mathbb{N}$-transitive sequence based on $n_i$ copies of $\mathbb{T}^{i}$  and
on $l_j$ copies of $\mathbb{T}^{\kappa_j}$ (where each $\kappa_j$ is an infinite cardinal with 
 $\kappa_j \le \mathfrak{c}$) for all $i,j$ taken simultaneously.
\end{ex}

\begin{rem}\label{zien}
Let $\mathbb{P} \subset \mathbb{N}$, and for each $n \in \mathbb{P}$, let $Z_n$ be a separable space. Let 
$Z:= \prod_{n \in \mathbb{P}} Z_n$. We know by the proof of Theorem~\ref{tabaclaro-v59} that
 $(Z^{\mathbb{Z}}, \Sigma, (w_n))$ is totally transitive for a certain point $(w_n) \in Z^{\mathbb{Z}}$ (where 
$\Sigma$ is the bilateral shift map). Clearly, if $u = (u_n) \in Z^{\mathbb{Z}}$, then  for every $n \in \mathbb{Z}$, $u_n$ can 
be written as $u_n= \pl u_n (m) \pr_{m \in \mathbb{P}}$, where $u_n (m) \in Z_m$. It is easy to check that 
$(\Sigma u)_n (m) = u_{n+1} (m)$ for every $n \in \mathbb{Z}$ and $m \in \mathbb{P}$. Consequently, if for 
$m \in \mathbb{P}$,  $\Sigma_m $ is the bilateral shift map on $Z_m^{\mathbb{Z}}$, then it turns out that 
$(Z_m^{\mathbb{Z}}, \Sigma_m, (w_n (m))_{n \in \mathbb{Z}} )$ is  totally transitive. It is not difficult to see that 
$(Z_m^{\mathbb{Z}}, \Sigma_m, (w_n (m))_{n \in \mathbb{Z}} )_{m \in \mathbb{P}}$ is $\mathbb{N}$-transitive due to its relation with $(Z^{\mathbb{Z}}, \Sigma, (w_n))$.
\end{rem}

\begin{ex}\label{zentu}
Remark~\ref{zien} gives us a way to find an $\mathbb{N}$-transitive sequence based on an arbitrary (countable) number of copies of spaces 
such as the Cantor set $\mathbf{K}$,  $\mathbb{T}^{\mathbb{Z}}$, $\mathbb{T}^{\mathbb{R}}$,
$\mathbb{N}^\mathbb{Z}$, 
$\mathbb{N}^\mathbb{R}$, $\ell^2$ (which is homeomorphic to $\mathbb{R}^\mathbb{Z}$), and any other
constructed in a similar way, that is,  $Z^{\kappa}$ (being homeomorphic to $\pl Z^{\kappa} \pr^{\mathbb{Z}}$) where $Z$ is separable and $\kappa$ is an infinite cardinal
less than or equal to  $\mathfrak{c}$.
\end{ex}

\begin{rem}\label{dindon}
Let $(Z, h, 1)$ be $\mathbb{L}$-transitive for $\mathbb{L} \subset \mathbb{N}$. Suppose that $Z$ is (homeomorphic to a)
dense subset of $Z_n$, and that $h$ can be extended continuously to a homeomorphism $h_n$ from $Z_n $ onto itself, 
for every $n \in \mathbb{P} \subset \mathbb{N}$. Then the sequence $ (Z_n, h_n, 1_n)_{n \in \mathbb{P}}$ is $\mathbb{L}$-transitive.  

We can state this fact more in general. For this, 
we say that given two semiflows  $(Z, h)$, $(Y, k)$, and two points $1 \in Z$, $1' \in Y$,  $(Z, h, 1)$ is a {\em dense restriction} of $(Y, k, 1')$ if  $Z$ is a dense subset of $Y$, $1=1'$, 
and $h$ is the restriction of $k$ to $Z$.

 Consider $\mathbb{P}_1, \mathbb{P}_2 \subset \mathbb{N}$. 
For each $n \in \mathbb{P}_1$ and $m \in \mathbb{P}_2$, let $(Z_n, h_n)$ and $(Y_m, k_m)$ be  semiflows, and let
 $1_n \in Z_n$  for every $n \in \mathbb{P}_1$. Suppose that  
 for every $n \in \mathbb{P}_1$,  there exists $m \in \mathbb{P}_2$ such that $(Z_n, h_n, 1_n)$ is a dense restriction 
of  $(Y_m, k_m, 1_{m_n})$. On the other hand,  suppose that the converse also holds, that is, for each $m \in \mathbb{P}_2$, there exists 
$n \in \mathbb{P}_1$ such that $(Z_n, h_n, 1_n)$ is a dense restriction of  $(Y_m, k_m, 1_{m_n})$.  Then, for $\mathbb{L} \subset \mathbb{N}$, we have 
that $(Z_n, h_n, 1_n)_{n \in \mathbb{P}_1 }$
 is $\mathbb{L}$-transitive if and only if 
$ (Y_m, k_m, 1_{m_n})_{m \in \mathbb{P}_2 }$ is.
\end{rem}

\begin{ex}
Taking into account Remark~\ref{dindon}, we see that if in Examples~\ref{zarrau} and \ref{zentu}, we define
 $Z_n^1=  \mathrm{Orb}^+ (h_n, 1_n)$ (or any subset containing this semiorbit), then both  $(Z_n^1, h_n, 1_n)_{n \in \mathbb{N} }$ and 
$(\beta Z_n^1, h_n, 1_n)_{n \in \mathbb{N} }$ are
$\mathbb{N}$-transitive.
\end{ex}

\section{$X$ needs not to be separable}\label{lucia10-m-07-rememsoper}

 Not much is known about the possibility of finding a nonseparable space $X$ such that $C(X)$ admits an isometric 
shift since the problem was proposed. Interesting results in this direction 
say that such an $X$ must have the countable chain condition (see Remark~\ref{dia6-4-07}), or even that   
$C_0 (X \setminus {\rm cl}_{X}(\mathcal{N}))$  (the space of $\mathbb{K}$-valued continuous functions vanishing at infinity) must have 
cardinality at most equal to $\mathfrak{c}$ (see \cite[Theorem 1.9]{GN}).

In this section, we prove that such an $X$ exists. As stated in Theorem~\ref{podmo-skilo-11m}, the space $X$ will be the topological sum of the maximal ideal space $\mathfrak{M}$ of $L^{\infty} (\mathbb{T})$ and $\mathcal{N} \cup \{\pmb{\infty}\}$. It is well known that $\mathfrak{M}$ is extremally disconnected, that is, the closure of each open subset is also open. In fact, each measurable subset $A$ of $\mathbb{T}$ determines via the Gelfand transform an open  and closed subset $\mathbf{G} (A)$ of $\mathfrak{M}$, and the sets obtained in this way form a basis for its topology (see \cite[p. 170]{H}). Now it is straightforward  to see that $\mathfrak{M}$ is not separable:  Let $(x_n)$ be a sequence in $\mathfrak{M}$, and consider
a partition (a.e) of $\mathbb{T}$ by $k$ arcs of equal length, $k \ge 3$. This determines a partition of $\mathfrak{M}$ into $k$ closed and open subsets of $\mathbb{T}$, and we select the arc $A_1$ such that $\mathbf{G} (A_1)$ contains $x_1$. Next do the same process with $k^2$ arcs of equal length, and pick $A_2$ with  $x_2 \in \mathbf{G} (A_2)$. Repeat the process infinitely many times, in such a way that each time we take $A_n$ of length $1/k^n$ such that $x_n \in \mathbf{G} (A_n)$. It is clear that if $A:= \bigcup_{n=1}^{\infty} A_n$, then $m \pl A \pr < 2 \pi$, so $\mathbf{G} (\mathbb{T} \setminus A)$ is a nonempty closed and open subset of $\mathfrak{M}$ containing no point $x_n$.

Notice that, since $\mathfrak{M}$ is not separable, any isometric shift on $C(\mathfrak{M})$ must be  of type $\mathrm{I}_0$. But there are none because  $\mathfrak{M}$ has no isolated points. Even more, in \cite[Corollary 2.5]{GHJR}, it is proved that no space $L^{\infty} (Z, \Sigma, \mu)$ admits an isometric shift if $\mu$ is non-atomic.

 As usual, we consider $\mathbb{T}$ oriented counterclockwise, and denote by $A(\alpha, \beta)$ the (open) arc of $\mathbb{T}$ beginning at $e^{i \alpha}$ and ending at $e^{i  \beta}$.

\begin{proof}[Proof of Theorem~\ref{podmo-skilo-11m}]
We start by defining a linear and surjective isometry on $L^{\infty} (\mathbb{T})$. We first consider the rotation $\psi (z) := z e^i$ for every $z \in \mathbb{T}$, and then define the  isometry $S: L^{\infty} (\mathbb{T}) \ra L^{\infty} (\mathbb{T})$ as $Sf:= - f \circ \psi$ for every $f \in L^{\infty} (\mathbb{T})$. On the other hand, 
using the Gelfand transform we have that the Banach algebra $L^{\infty} (\mathbb{T})$  is  isometrically isomorphic to $C(\mathfrak{M})$, so $S$
 determines a linear and surjective isometry $T_S : C(\mathfrak{M}) \ra C(\mathfrak{M})$. Also, by the Banach-Stone theorem, there exists a homeomorphism $\phi: \mathfrak{M} \ra \mathfrak{M}$ such that $T_S f = - f \circ \phi$ for every $f \in C(\mathfrak{M})$. Notice that this is valid both in the real and complex cases (see for instance \cite[p. 187]{G}).

Let $X:= \mathfrak{M} + \mathcal{N} \cup \{\pmb{\infty}\}$. The definition of $T_S$ can 
be extended to a new isometry  $T: C(X) \ra C(X)$ in three steps.  First, for each $f \in C(X)$, we put $(T f) (x) := (T_S f) (x)$ if 
$x \in \mathfrak{M}$. In the same way $(Tf) (\mathbf{n}) := (f \circ \phi) (\mathbf{n})$ if $\mathbf{n} \in \mathcal{N} \cup \{\pmb{\infty}\} \setminus \{\mathbf{1}\}$ (where $\phi: \mathcal{N} \setminus \{\mathbf{1}\} \ra \mathcal{N}$ is the canonical map 
sending each $\mathbf{n}$ into $\mathbf{n-1}$, which obviously can be extended as   
$\phi (\pmb{\infty}) := \pmb{\infty}$). Finally, we put $$(Tf) (\mathbf{1} ):= \frac{1}{2 \pi}\int_{A(0, 2 \pi \Phi)} f dm,$$ where   $\Phi:= \pl \sqrt{5}-1 \pr /2$ is the golden ratio conjugate. 
It is easy to verify that $T$ is a codimension one linear isometry, so we just need to prove that $\bigcap_{i=1}^{\infty} T^i (C(X)) =\{\mathbf{0}\}$.   

Suppose then that  $f \in \bigcap_{i=1}^{\infty} T^i (C(X))$. It is  easy to check that 
\begin{eqnarray*}
f(\mathbf{n}) &=&  (T^{-n +1} f) (\mathbf{1}) \\
&=&  (T(T^{-n} f)) (\mathbf{1}) \\
&=& \frac{(-1)^{n}}{2 \pi} \int_{A(0, 2 \pi \Phi)} f \circ \psi^{-n} dm \\
&=& \frac{(-1)^{n}}{2 \pi} \int_{A(n, n + 2 \pi \Phi )} f  dm.
\end{eqnarray*}

On the other hand, if we fix any $\alpha \in \mathbb{T}$, then there exist two increasing sequences $(n_k)$ and $(m_k)$ in $2 \mathbb{N}$ and $2 \mathbb{N} +1$, respectively, 
converging to $\alpha   \mod 2 \pi$. An easy application of the  Dominated Convergence Theorem proves that 
$ \int_{A(\alpha, \alpha + 2 \pi \Phi)} f  dm = 2 \pi \lim f (\mathbf{n}_k) $, and $\int_{A(\alpha, \alpha + 2 \pi \Phi )} f  dm  = - 2 \pi \lim f (\mathbf{m}_k)$. By continuity, we 
deduce that $$\int_{A(\alpha, \alpha + 2 \pi \Phi )} f  dm = 2 \pi f(\pmb{\infty}) = - \int_{A(\alpha, \alpha + 2 \pi \Phi )} f  dm.$$ Obviously, this implies that $\int_{A(\alpha, \alpha + 2 \pi \Phi)} f  dm = 0$ for every $\alpha \in \mathbb{T}$, 
 and $f(\pmb{\infty}) =0$. In particular this proves that $f (\mathbf{n}) = 0$ for every 
$\mathbf{n} \in \mathcal{N}$. As a consequence we can identify  $f \in \bigcap_{n=1}^{\infty} T^n (C(X))$ with an element  $f \in L^{\infty} (\mathbb{T})$ satisfying $\int_{A(\alpha, \alpha + 2 \pi \Phi )} f  dm= 0$ for every $\alpha \in \mathbb{T}$. On the other hand, it is clear that we may assume that $f$ takes values just in $\mathbb{R}$.

\medskip

{\sc Claim.} {\em $\int_{A(\alpha, \alpha + 2 \pi \Phi^n)} f dm = (-1)^n F(n-1) \int_{\mathbb{T}} f dm$ for every $\alpha \in \mathbb{T}$ and $n \in \mathbb{N}$, where $F(n)$ denotes the $n$-th Fibonacci number.}

Let us prove the claim inductively on $n$. We know that it holds for $n=1$. Also notice that $\Phi+ \Phi^2= 1$, so $\Phi^n + \Phi^{n+1} = \Phi^{n-1}$ for every $n \in \mathbb{N}$. 

The case $n=2$ is immediate because, since $$\mathbb{T} = A \pl \alpha, \alpha+  2 \pi \Phi^2 \pr  \cup A \pl \alpha+  2 \pi \Phi^2, \alpha+  2 \pi \pl \Phi +  \Phi^2 \pr \pr$$ a. e., then we have $\int_{\mathbb{T}} f  dm=  \int_{A \pl \alpha, \alpha+   2 \pi \Phi^2 \pr} f  dm$ for every $\alpha \in \mathbb{T}$. 

Now assume that, given $k \ge 2$, the claim is true for every $n \le k$. Then we see that, for any  $\alpha \in \mathbb{T}$, 
$$A \pl \alpha, \alpha+  2 \pi \Phi^{k-1} \pr = 
A \pl \alpha, \alpha+  2 \pi  \Phi^{k+1} \pr \cup A \pl \alpha+  2 \pi  \Phi^{k+1}, 
    \alpha+  2 \pi \pl \Phi^{k}  +   \Phi^{k+1} \pr \pr  $$a. e., so
$$ (-1)^{k-1} F(k-2) \int_{\mathbb{T}} f dm = \int_{A \pl \alpha, \alpha+  2 \pi  \Phi^{k+1} \pr} f dm + 
(-1)^{k} F(k-1) \int_{\mathbb{T}} f dm,$$
and the conclusion proves the claim.

\medskip

The claim, combined with the fact that $f$ is essentially bounded,  implies that $\int_{\mathbb{T}} f  dm =0$, and  consequently $\int_{A(\alpha, \alpha + 2 \pi \Phi^n)} f dm =0$ for every $\alpha \in \mathbb{T}$ and every $n \in \mathbb{N}$.

Now, it is easy to see that if $U$ is an open subset of $\mathbb{T}$, then $U$ is the union of countably many pairwise disjoint arcs whose lengths belong to the set $\{2 \pi \Phi^n : n \in \mathbb{N}\}$. Now, applying again  the Dominated Convergence Theorem, we see that $\int_U f dm=0$. Obviously, this implies that $\int_K f dm =0$ whenever $K \subset \mathbb{T}$ is compact.

Finally take $C^+ := \{z \in \mathbb{T} : f (z) > 0\}$. We know that there exists a sequence of compact subsets $K_n$ of $C^+$, with $K_n \subset K_{n+1}$ for every $n \in \mathbb{N}$, and such that $\lim_{n \ra \infty} m( C^+ \setminus K_n) =0$. Clearly, the above fact and the Monotone Convergence Theorem imply that $\int_{C^+} f dm=0$, and then $m (C^+) =0$. Now we can easily conclude that $f \equiv 0 $ a. e., and consequently $T$ is a shift.
\end{proof}

Next we prove Theorem~\ref{sinregalocumple}. It provides  nonseparable examples with $2^{\mathfrak{c}}$  infinite 
connected  components, each homeomorphic to a (finite or infinite dimensional)  torus: To show it, we just need to see that the 
cardinal of $\mathfrak{M}$ is $2^{\mathfrak{c}}$. This follows from the fact that $\mathfrak{M}$ is homeomorphic to an 
infinite closed subset of $\beta \mathbb{N} \setminus \mathbb{N}$, and consequently its cardinal must be $2^{\mathfrak{c}}$ (see \cite{N} and \cite[Corollary 9.2]{GJ}).

\begin{proof}[Proof of Theorem~\ref{sinregalocumple}]
Write the isometric shift $T: C( \mathfrak{M} + \mathcal{N} \cup \{\pmb{\infty}\}) \ra C( \mathfrak{M} + \mathcal{N} \cup \{\pmb{\infty}\})$ given in the proof of Theorem~\ref{podmo-skilo-11m} as $T= T[a, \phi, \Delta] $. 
Obviously,  $\Delta \equiv 0$ on $C( \mathcal{N} \cup \{\pmb{\infty}\})$, and it can be considered as an element of $C( \mathfrak{M})'$.

 Consider a subset $\mathbb{P}$ of $\mathbb{R}$
 with cardinal equal to $\kappa$, and suppose that  $\{1/2 \pi\} \cup \{\rho_{\alpha} :\alpha \in \mathbb{P}\}$ is a family  of real numbers linearly independent over $\mathbb{Q}$. Then put
 $\rho_{\kappa} := [\rho_{\alpha}]_{\alpha \in \mathbb{P}}$  (see Example~\ref{dominicus}). 

Define $\phi_{\kappa} : \mathfrak{M} \times \mathbb{T}^{\kappa} \ra \mathfrak{M} \times \mathbb{T}^{\kappa}$ as
$$\phi_{\kappa} (x, \mathbf{z}) := \pl \phi (x), \rho_{\kappa} (\mathbf{z}) \pr$$ for every $x \in \mathfrak{M}$, and  
$\mathbf{z} \in \mathbb{T}^{\kappa}$. 

Select now  a point $\mathbf{v_{\kappa}}$   in 
$\mathbb{T}^{\mathbb{P}} = \mathbb{T}^{\kappa}$, and consider  the evaluation map 
$\Gamma_{\mathbf{v}_{\kappa}} \in C(\mathbb{T}^{\kappa})'$. Both  $\Delta$ and $\Gamma_{\mathbf{v}_{\kappa}}$ are {\em
 positive} linear functionals, and so is the product $\Delta \times \Gamma_{\mathbf{v}_{\kappa}} \in C( \mathfrak{M} \times \mathbb{T}^{\kappa})'$,
 which also satisfies $\vc \Delta \times \Gamma_{\mathbf{v}_{\kappa}} \vd \le 1$ (see \cite[\textsection 13]{HR} for details). 

Given $f \in C(\mathfrak{M} \times \mathbb{T}^{\kappa})$ and $\mathbf{z} \in \mathbb{T}^{\kappa}$, we write 
$f_{\mathbf{z}} :  \mathfrak{M}  \ra \mathbb{K}$ meaning $f_{\mathbf{z}} (x) := f(x, \mathbf{z})$ 
for every $x \in \mathfrak{M}$. Obviusly $f_{\mathbf{z}} $ belongs to  $C( \mathfrak{M})$, and 
$(\Delta \times \Gamma_{\mathbf{v}_{\kappa}}) (f) = \Gamma_{\mathbf{v}_{\kappa}} \pl \Delta (f_{\mathbf{z}}) \pr = \Delta (f_{\mathbf{v}_{\kappa}})$.

Now, for  $X_{\kappa} := \mathfrak{M} \times \mathbb{T}^{\kappa} + \mathcal{N} \cup \{\pmb{\infty}\}$,  define $a_{\kappa} \in C(X_{\kappa} \setminus \{\mathbf{1}\})$ as $a_{\kappa} \equiv -1$ on $\mathfrak{M} \times \mathbb{T}^{\kappa}$, and $a_{\kappa} \equiv 1$ everywhere 
else, and put $T_{\kappa} := T[a_{\kappa}, \phi_{\kappa}, \Delta \times \Gamma_{\mathbf{v}_{\kappa}}]$. 

Let $\psi: \mathbb{T} \ra \mathbb{T}$ and $\Phi$ be as in the proof of Theorem~\ref{podmo-skilo-11m}. Given  $f \in \bigcap_{i=1}^{\infty} T_{\kappa}^i (X_{\kappa})$,   we have that for every $k \in \mathbb{N}$, 
\begin{eqnarray*}
f(\mathbf{k}) &=&  (T_{\kappa}^{-k +1} f) (\mathbf{1}) \\
&=&  \pl \Delta \times \Gamma_{\mathbf{v}_{\kappa}} \pr \pl T_{\kappa}^{-k}  f \pr \\
&=& (-1)^{k} \pl \Delta \times \Gamma_{\mathbf{v}_{\kappa}} \pr \pl f \circ \phi_{\kappa}^{-k} \pr \\
&=& (-1)^{k}  \Delta \pl f \circ \phi_{\kappa}^{-k} \pr_{\mathbf{v}_{\kappa}} \\
&=& (-1)^{k}  \Delta \pl f_{\rho_{\kappa}^{-k} \pl \mathbf{v}_{\kappa} \pr} \circ \phi^{-k} \pr \\
&=& \frac{(-1)^{k}}{2 \pi} \int_{A(0, 2 \pi \Phi)} f_{\rho_{\kappa}^{-k}  (\mathbf{v}_{\kappa})} \circ \psi^{-k} dm \\
&=& \frac{(-1)^{k}}{2 \pi} \int_{A(k, k + 2 \pi \Phi )} f_{\rho_{\kappa}^{-k}  (\mathbf{v}_{\kappa})}  dm .
\end{eqnarray*}

To continue with the proof, we need an elementary result:

\medskip

{\sc Claim.} {\em 
Suppose that $(\mathbf{z}_{\lambda})_{\lambda \in D}$ is a net in $\mathbb{T}^{\kappa}$ converging to $\mathbf{z}_0$. Then $\lim_{\lambda} \vc f_{\mathbf{z}_{\lambda}} - f_{\mathbf{z}_0}  \vd =0$. }

Let us prove the claim. If it is not true, then there is an $\epsilon>0$ such that, for every $\lambda \in D$, there exists $\nu \in D$, $\nu \ge \lambda$,  such that $\vc f_{\mathbf{z}_{\nu}} - f_{\mathbf{z}_0}  \vd \ge \epsilon$. It is easy to see that the set $E$ of all $\nu \in D$ satisfying the above inequality is a directed set, and that
 $(\mathbf{z}_{\nu})_{\nu \in E}$  is a subnet of $(\mathbf{z}_{\lambda})_{\lambda \in D}$. Moreover there is a net  
 $(x_{\nu})_{\nu \in E}$ in $\mathfrak{M}$
such that $\va f(x_{\nu}, \mathbf{z}_{\nu}) - f(x_{\nu}, \mathbf{z}_0) \vb \ge \epsilon$ for every $\nu \in E$. Since 
$\mathfrak{M} \times \mathbb{T}^{\kappa}$ is compact, there exist a point $(x_0 , \mathbf{z}_0') \in \mathfrak{M} \times \mathbb{T}^{\kappa}$ and a subnet $(x_{\eta}, \mathbf{z}_{\eta})_{\eta \in F}$ of 
$(x_{\nu}, \mathbf{z}_{\nu} )_{\nu \in E}$ converging to $(x_0, \mathbf{z}_0' ) $. Obviously $(\mathbf{z}_{\eta})_{\eta \in F}$ is a subnet of $(\mathbf{z}_{\nu})_{\nu \in E}$, so  $\mathbf{z}_0 =\mathbf{z}_0'$. Consequently both $\pl x_{\eta}, \mathbf{z}_{\eta} \pr_{\eta \in F}$ and $\pl x_{\eta}, \mathbf{z}_0 \pr_{\eta \in F}$ converge to $\pl x_0  , \mathbf{z}_0 \pr$.  Taking limits, this 
implies  $\va f (x_0, \mathbf{z}_0) - f (x_0, \mathbf{z}_0)\vb \ge \epsilon$, which is absurd.

\medskip

Now, fix $(\alpha, \mathbf{w}) \in \mathbb{T} \times \mathbb{T}^{\kappa}$ and $\epsilon >0$. We know that 
$(\alpha, \mathbf{w})$ belongs to the closure of both $$\mathbf{N}_j:= \tl \pl e^{i n}, \rho_{\kappa}^{-n} 
(\mathbf{v}_{\kappa}) \pr :n \in 2 \mathbb{N} - j \tr \setminus \tl (\alpha, \mathbf{w})\tr ,$$ $j=0,1$ (see Example~\ref{dominicus}). We first consider the case 
$j=0$, and take a net $(y_{\lambda})_{\lambda \in D} = (e^{i n_{\lambda}}, \rho_{\kappa}^{-{n_{\lambda}}} (\mathbf{v}_{\kappa}))_{\lambda \in D}$ in $\mathbf{N}_0$ converging to $(\alpha, \mathbf{w})$. Since $\pl e^{i n_{\lambda}} \pr_{\lambda \in D}$ converges to $\alpha$, there exists $\lambda_1\in D$ such that $$\va \int_{A(\alpha, \alpha + 2 \pi \Phi )} f_{\mathbf{w}}  dm  - \int_{A(n_{\lambda}, n_{\lambda} + 2 \pi \Phi )} f_{\mathbf{w}} \vb < \frac{\epsilon}{2}$$ for every $\lambda \ge \lambda_1$.

On the other hand, by the claim, there exists $\lambda_2 \in D$ such that, if $\lambda \ge \lambda_2$, then 
$\vc  f_{\mathbf{w}} - f_{\rho_{\kappa}^{-n_{\lambda}}  (\mathbf{v}_{\kappa})}   \vd < \epsilon /4 \pi$, so  
 $$\va \int_{A(n_{\nu} , n_{\nu} + 2 \pi \Phi )} f_{\mathbf{w}}  dm  - \int_{A(n_{\nu}, n_{\nu} + 2 \pi \Phi )} f_{\rho_{\kappa}^{-n_{\lambda}}  (\mathbf{v}_{\kappa})}  dm   \vb < \frac{\epsilon}{2}$$ for every $\nu \in D$.   
We easily deduce that 
$$\lim_{\lambda} \int_{A(n_{\lambda}, n_{\lambda} + 2 \pi \Phi )} f_{\rho_{\kappa}^{-n_{\lambda}}  (\mathbf{v}_{\kappa})}  dm = \int_{A(\alpha, \alpha + 2 \pi \Phi )} f_{\mathbf{w}}  dm, $$
 and consequently $2 \pi f(\pmb{\infty}) = \int_{A(\alpha, \alpha + 2 \pi \Phi )} f_{\mathbf{w}}  dm$. In a similar way,  working with $\mathbf{N}_1$,  we see that $2 \pi  f(\pmb{\infty}) = - \int_{A(\alpha, \alpha + 2 \pi \Phi )} f_{\mathbf{w}}  dm$. With  the same arguments as in the proof of Theorem~\ref{podmo-skilo-11m}, we conclude that $f_{\mathbf{w}} \equiv 0$, and finally $f \equiv 0$, as we wanted to prove.
\end{proof}

\begin{proof}[Proof of Theorem~\ref{notatillo2}]
Notice first that $L^{\infty} (\mathbb{T})$ is isometrically isomorphic to $L^{\infty} (\mathbb{T}_1 \cup \mathbb{T}_2)$,
where $\mathbb{T}_i$, $i=1,2$, are disjoint copies of $\mathbb{T}$ endowed with the Lebesgue measure. It is not hard to see that this implies that $C( \mathfrak{M})$ and $C( \mathfrak{M} + \mathfrak{M})$ are isometrically isomorphic, so $\mathfrak{M}$ and $\mathfrak{M} + \mathfrak{M}$ are homeomorphic. Assume that   $T= T[a, \phi, \Delta]$ is  the isometric shift given in the proof of Theorem~\ref{podmo-skilo-11m}. We first define a homeomorphism $\chi : \mathfrak{M} \times \{0,1\} \ra \mathfrak{M} \times \{0,1\}$ as $\chi (x, i) = (\phi (x), i+1 \mod 2)$ for every $(x,i)$. 
For $i=0,1$, and $f \in C(\mathfrak{M} \times \{0,1\})$, denote by  $f \times \{i\} $ its restriction to $\mathfrak{M} \times \{i\}$, and put $\Delta_i (f) := \Delta (f \times \{i\})$.

Let $\rho_{\kappa} : \mathbb{T}^{\kappa} \ra \mathbb{T}^{\kappa} $,  $\mathbf{v_{\kappa}}$, and $\Gamma_{\mathbf{v_{\kappa}}}$  be as in the proof of Theorem~\ref{sinregalocumple}.

Finally consider $X_{\kappa}:= \mathfrak{M} \times \{0,1\} + \mathbb{T}^{\kappa} + \mathcal{N} \cup \{\pmb{\infty}\}$, and define $T_{\kappa}: C(X_{\kappa}) \ra
C(X_{\kappa}) $ to be $T_{\kappa} := T [a_{\kappa} , \phi_{\kappa}, \Delta_{\kappa}]$, where
\begin{itemize}
\item $a_{\kappa} \equiv -1$ on $\mathfrak{M} \times \{0\} \cup \mathbb{T}^{\kappa}$,  and $a_{\kappa} \equiv 1$ everywhere else.
\item $\phi_{\kappa} = \chi$ on $\mathfrak{M} \times \{0,1\}$, and $\phi_{\kappa} = \rho_{\kappa}$ on $\mathbb{T}^{\kappa}$.
\item $\Delta_{\kappa}  := \pl \Delta_0  + \Delta_1 + \Gamma_{\mathbf{v}_{\kappa}} \pr / 3$.
\end{itemize}

As above, if $f \in \bigcap_{n=1}^{\infty} T_{\kappa}^n (C(X_{\kappa}))$,  $k \in \mathbb{N}$, and $$\tau (k) := \frac{k(k-1) \mod 4}{2} ,$$ then 

\begin{eqnarray*}
3 f(\mathbf{k}) &=&  3 \pl T_{\kappa}^{-k +1} f \pr \pl \mathbf{1} \pr \\
&=&  \Delta_0 \pl T_{\kappa}^{-k} f \pr + \Delta_1 \pl T_{\kappa}^{-k} f \pr  + \Gamma_{\mathbf{v}_{\kappa}} \pl T_{\kappa}^{-k} f \pr\\
&=&  \Delta \pl \pl T_{\kappa}^{-k} f \pr  \times \{0\}\pr  + \Delta \pl \pl  T_{\kappa}^{-k} f\pr  \times \{1\}\pr   +  \pl T_{\kappa}^{-k} f \pr \pl  \mathbf{v}_{\kappa} \pr \\
&=& (-1)^{\tau (k)}  \Delta \pl \pl  f    \times \{k \mod 2\} \pr \circ \phi^{-k} \pr \\ 
&&+ (-1)^{\tau (k+1)} \Delta \pl \pl   f   \times \{k + 1 \mod 2\} \pr \circ \phi^{-k}  \pr  \\ 
&& +  (-1)^k \pl f \circ \rho_{\kappa}^{-k}  \pr \pl  \mathbf{v}_{\kappa} \pr \\
&=& \frac{(-1)^{\tau (k)}}{2 \pi} \int_{A(0, 2 \pi \Phi)} \pl f \times \{k \mod 2\} \pr  \circ \psi^{-k} dm \\
&& + \frac{(-1)^{\tau (k+1)}}{2 \pi} \int_{A(0, 2 \pi \Phi)} \pl f \times \{k +1 \mod 2\} \pr  \circ \psi^{-k} dm \\
&&+ (-1)^k  f \pl  \rho_{\kappa}^{-k}  \pl  \mathbf{v}_{\kappa} \pr \pr \\
&=& \frac{(-1)^{\tau (k)}}{2 \pi} \int_{A(k, k + 2 \pi \Phi)}  f \times \{k \mod 2\}  dm \\
&&+ \frac{(-1)^{\tau (k+1)}}{2 \pi} \int_{A(k, k+ 2 \pi \Phi)}  f \times \{k +1 \mod 2\}  dm \\
&&+ (-1)^k  f \pl  \rho_{\kappa}^{-k}  \pl  \mathbf{v}_{\kappa} \pr \pr . 
\end{eqnarray*}

Next fix  $\alpha \in \mathbb{T}$, $\mathbf{w}  \in \mathbb{T}^{\kappa}$, and for $j=0, 1, 2, 3$, take increasing sequences $(n_k^j)$ in $4 \mathbb{N} + j$ such that 
$\lim_{k \ra \infty} n_k^j = \alpha \mod 2 \pi$, and $\lim_{k \ra \infty} f \pl \rho_{\kappa}^{-n_k^j}  \pl  \mathbf{v}_{\kappa} \pr \pr = f \pl \mathbf{w} \pr$.
 Now put $$\mathbf{X}_i^{\alpha} := \frac{1}{2 \pi} \int_{A(\alpha , \alpha  + 2 \pi \Phi)}  f \times \{i\}  dm$$ for $i=0, 1$. 
Taking into account that $\tau (n_k^j)$ is constant for each $j$, and that $\tau (2) = 1 = \tau(3)$, and 
$\tau(1) = 0 = \tau (4)$, we have that the following equalities hold:

$$\begin{array}{lcrcl}
3 f(\mathbf{\pmb{\infty}}) &=&   \mathbf{X}_1^{\alpha} -\mathbf{X}_0^{\alpha}-  f \pl \mathbf{w}  \pr && (\mathrm{case } \hspace{.03in} j= 1) \\
&=& - \mathbf{X}_0^{\alpha} -\mathbf{X}_1^{\alpha} + f \pl \mathbf{w}  \pr  &&  ( \mathrm{case } \hspace{.03in} j= 2)   \\
&=&  - \mathbf{X}_1^{\alpha} + \mathbf{X}_0^{\alpha}-  f \pl \mathbf{w}   \pr  && ( \mathrm{case } \hspace{.03in} j= 3)  \\
&=&  \mathbf{X}_0^{\alpha} + \mathbf{X}_1^{\alpha} + f \pl \mathbf{w} \pr && ( \mathrm{case } \hspace{.03in} j= 0)     . 
\end{array}
$$

We  deduce that $\mathbf{X}_i^{\alpha} =0$ for every $\alpha \in \mathbb{T}$ and $i =0,1$, and that $f \equiv 0$ on $\mathbb{T}^{\kappa}$. As in the proof 
of Theorem~\ref{podmo-skilo-11m}, we easily conclude that $f \equiv 0$.  
\end{proof}

\section{Complex $n$-generated shifts and nonseparable spaces with a finite number of infinite connected components}\label{emphado}

This  is the only section of the paper where $\mathbb{K}$ is always taken equal to $\mathbb{C}$. Theorem~\ref{notatillo} 
tells us that in some cases  we can obtain an $(n+1)$-generated shift (obviously  of type $\mathrm{I}_0$) from an $n$-generated one just by adding a connected component 
($n \in \mathbb{N}$). A similar proof shows that it is possible to obtain nonseparable examples with arbitrary (finitely)  many  infinite connected components. For the different behavior in the real setting, see Example~\ref{torero}.

\begin{thm}\label{notatillo}
Let $\mathbb{K} = \mathbb{C}$.  Suppose that $n \ge 2$, and that $(\kappa_i)_{i=1}^{n}$ is a finite sequence of cardinals
satisfying $0 \le \kappa_i \le \mathfrak{c}$ for every $i$. 
Then  there exists  an $n$-generated isometric shift $T_n: C_{\mathbb{C}} (X_n) \ra C_{\mathbb{C}} (X_n)$, where $X_n = \mathbb{T}^{\kappa_1} + \ldots + \mathbb{T}^{\kappa_n} + \mathcal{N} \cup \{\pmb{\infty}\}$. 
\end{thm}

\begin{proof}
We will prove it inductively on $n$. We start at $n=1$ (which is a special case). Let $\mathbb{P}_1, \ldots, \mathbb{P}_n, \ldots$ be pairwise disjoint subsets of 
$\mathbb{R}$  of cardinalities $\kappa_1, \ldots, \kappa_n, \ldots $, respectively. Consider then a family 
$\Lambda:= \{\rho_{\alpha} : \alpha \in \mathbb{R}\}$ of real numbers linearly independent over $\mathbb{Q}$, and put
 $\sigma_i := [\rho_{\alpha}]_{\alpha \in \mathbb{P}_i}$ for each $i \le n$. Also let $\mathbf{v}_i $ be a point in 
$\mathbb{T}^{\mathbb{P}_i}$ (see Example~\ref{dominicus}).   

Next suppose that  $X_n := \mathbb{T}^{\mathbb{P}_n} + \ldots + \mathbb{T}^{\mathbb{P}_1} + \mathcal{N} \cup \{\pmb{\infty}\}$, and define  $\phi_n : X_n \ra X_n$ as $\sigma_i$ on each $\mathbb{T}^{\mathbb{P}_i}$. For  $n \in \mathbb{N}$, let $z_n \in \mathbb{C}\setminus \{0\}$, with $\va z_n \vb \le 1/2^n$, and    $\zeta_n := e^{i \pi/2^{n-1}}$. 
 Define a codimension $1$ linear isometry $T_n$ on $C_{\mathbb{C}} (X_n)$ as $T_n := T [a_n, \phi_n, \Delta_n]$, where $a_n \equiv \zeta_i$
 on $\mathbb{T}^{\mathbb{P}_i}$  for each $i \le n$, and $a_n \equiv 1$ on $ \mathcal{N} \cup \{\pmb{\infty}\}$, and where $\Delta_n (f) := \sum_{i=1}^{n} z_i f(\mathbf{v}_i)$ for every $f$.

It is easy to see that $T_1 : C_{\mathbb{C}} (X_1) \ra C_{\mathbb{C}} (X_1)$ is an isometric shift (both of type $\mathrm{I}$ and type $\mathrm{II}$). Now let us show that if $T_i$ are $i$-generated isometric shifts for $i \le n$, then $T_{n+1}$ is  an $(n+1)$-generated isometric shift.

Suppose  that  $f \in \bigcap_{m=1}^{\infty} T_{n+1}^m (C_{\mathbb{C}} (X_{n+1}))$. It is  easy to check that 
\begin{eqnarray*}
f(\mathbf{k}) &=&  (T_{n+1}^{-k +1} f) (\mathbf{1}) \\
&=&  \sum_{i=1}^{n+1} z_i \pl T_{n+1}^{-k} f \pr(\mathbf{v}_i) \\
&=&  z_{n+1} \zeta_{n+1}^{-k} \pl f \circ \sigma_{n+1}^{-k}  \pr (\mathbf{v}_{n+1}) +  \sum_{i=1}^{n} z_i  \zeta_i^{-k} \pl f \circ \sigma_i^{-k}  \pr (\mathbf{v}_i) .
\end{eqnarray*}
whenever $k \in \mathbb{N}$.

Fix $x_1 \in \mathbb{T}^{\mathbb{P}_1}, \ldots, x_{n+1} \in \mathbb{T}^{\mathbb{P}_{n+1}}$. For $j=0, 1$, we can take increasing sequences $(n_k^j)$ in $2^{n+1} \mathbb{N}$ and $2^{n+1} \mathbb{N} + 2^{n}$, respectively, such that the sequences
$$\pl \pl f \circ \sigma_1^{-n_k^j}  \pr \pl \mathbf{v}_1 \pr, \ldots , \pl f \circ \sigma_{n+1}^{-n_k^j}  \pr (\mathbf{v}_{n+1}) \pr_{k \in \mathbb{N}}$$
converge to $(f(x_1), \ldots, f(x_{n+1})) \in \mathbb{C}^{n+1}$
, for $j=0, 1$.

This means, on the one hand, that 
\begin{eqnarray*}
f (\pmb{\infty}) &=& \lim_{k \ra \infty} f (\mathbf{n}_k^0) \\
&=& \lim_{k \ra \infty} z_{n+1} \zeta_{n+1}^{-n_k^0} \pl f \circ \sigma_{n+1}^{-n_k^0}  \pr (\mathbf{v}_{n+1}) +  \sum_{i=1}^{n} z_i  \zeta_i^{-n_k^0} \pl f \circ \sigma_i^{-n_k^0}  \pr (\mathbf{v}_i) \\
&=& z_{n+1} f (x_{n+1}) + \sum_{i=1}^{n} z_{i} f (x_{i}).
\end{eqnarray*}
And, on the other hand, 
\begin{eqnarray*}
f (\pmb{\infty}) &=& \lim_{k \ra \infty} f (\mathbf{n}_k^1) \\
&=& - z_{n+1} f (x_{n+1}) + \sum_{i=1}^{n} z_{i} f (x_{i}).
\end{eqnarray*}

We deduce that $f(x_{n+1}) =0$, that is, $f \equiv 0$ on $\mathbb{T}^{\mathbb{P}_{n+1}}$, and
consequently $f \in \bigcap_{m=1}^{\infty} T_n^m (C_{\mathbb{C}} (X_n))$. Since $T_n$ is a shift, we conclude that $f \equiv 0$ on $X_{n+1}$. It is also easy to see that $T_{n+1}$ is $(n+1)$-generated.
\end{proof}

\begin{thm}\label{jarboresquil}
Let $\mathbb{K} = \mathbb{C}$.  Suppose that $n \ge 1$, and that $(\kappa_i)_{i=1}^{n}$ is a finite sequence of cardinals
satisfying  $0 \le \kappa_i \le \mathfrak{c}$ for every $i $. Then there exists an isometric shift 
$T^{\mathfrak{M}}_n:  C_{\mathbb{C}} (X_n^{\mathfrak{M}}) \ra C_{\mathbb{C}} (X_n^{\mathfrak{M}})$, where $X_n^{\mathfrak{M}} = \mathfrak{M} + \mathbb{T}^{\kappa_1} + \ldots + \mathbb{T}^{\kappa_n} + \mathcal{N} \cup \{\pmb{\infty}\}$. 
\end{thm}

\begin{proof}
The proof is similar to that of Theorem~\ref{notatillo}. We consider the homeomorphism $\phi$ on $\mathfrak{M}$ coming from the rotation  $\psi: \mathbb{T} \ra \mathbb{T}$ given in  the proof of Theorem~\ref{podmo-skilo-11m}.

Fix $n \in \mathbb{N}$, and assume that $z_{n+1}$, $\zeta_{n+1}$, $X_n$, and $T_n =T [a_n, \phi_n, \Delta_n]$ are as  in the proof of Theorem~\ref{notatillo}.  We are going to define an isometric shift on $X_n^{\mathfrak{M}}$. First put  $$\Delta_n^{\mathfrak{M}}  := \frac{z_{n+1}}{2 \pi} \int_{A(0, 2 \pi \Phi )} f  dm + \Delta_n (f) .$$
Obviously we are assuming that $1/2 \pi$ does not belong to the linear span (over $\mathbb{Q}$) of $\{\rho_{\alpha} : \alpha \in \mathbb{P}_1 \cup \cdots \cup \mathbb{P}_n\}$. Let $a_{n}^{\mathfrak{M}} \in C
_{\mathbb{C}} (X_n^{\mathfrak{M}})$ be equal to $\zeta_{n+1}$ on $\mathfrak{M}$, and equal to $a_n$ on $X_n$, and let $\phi_n^{\mathfrak{M}} : X_n^{\mathfrak{M}} \ra X_n^{\mathfrak{M}}$ be defined as $\phi_n$ on $X_n$, and as $\phi$ on $\mathfrak{M}$.

 We consider  $T_{n}^{\mathfrak{M}} := T [a_{n}^{\mathfrak{M}}, \phi_n^{\mathfrak{M}} , \Delta_n^{\mathfrak{M}} ]$. 
Following the same process as in the proof of Theorem~\ref{notatillo}, we easily obtain that $0 = z_{n+1} \int_{A(\alpha, 2 \pi \Phi + \alpha)} f  dm$ for every $\alpha \in \mathbb{T}$. As in the proof of Theorem~\ref{podmo-skilo-11m}, we see that $f \equiv 0$ on $\mathbb{T}$, which is to say that $f \equiv 0$ on $\mathfrak{M}$. We deduce that $f \in \bigcap_{m=1}^{\infty} T_n^{m} (C_{\mathbb{C}} (X_n))$, and consequently $f \equiv 0$. 
\end{proof}

\begin{rem}\label{tillomalo}
Notice that in both Theorems~\ref{notatillo} and \ref{jarboresquil}, we allow the possibility that $\kappa_i = \kappa_j$ for some (or all) $i \neq j$.
\end{rem}

\section{A   method for obtaining  examples of $n$-generated shifts}\label{majorca}

In Section~\ref{emphado}, we saw in particular how to obtain inductively $n$-generated shifts in the complex setting, for $n \in \mathbb{N}$, $n \ge 2$. In this section, we will give a  method for constructing examples of  isometric shifts  of type $\mathrm{I}_0$ which are $n$-generated, allowing the case $n= \infty$. These operators will be based on separable spaces  consisting of the topological sum of a disconnected space and $\mathcal{N} \cup \{ \pmb{\infty}\}$.
In the case when $n= \infty$, we  construct these spaces as compactifications of topological sums with infinitely many  summands. Here we study simultaneously the real and complex settings (Example~\ref{torero} is the only result given just in the real context).

Example~\ref{torero} shows in fact that the procedure followed in Section~\ref{emphado} is no longer valid when dealing with $\mathbb{K} = \mathbb{R}$.

\begin{ex}\label{torero}
Let $\mathbb{K} = \mathbb{R}$. Suppose that $X= Y +  X_1 + X_2 + X_3$ is compact,
where each $X_i$ is connected and nonempty, and  $ \mathcal{N}  \subset Y$. 
Let $T = T[a, \phi, \Delta]$ be a codimension $1$ linear isometry  on $C_{\mathbb{R}} (X)$, and assume that
$\phi (X_i) = X_i$, $i=1,2,3$. Let us see that $T$ is not a shift. First, there are $i, j$, $i \neq j$,
with $a(X_i) = a \pl X_j \pr \in \{ -1, 1\}$. There are also $\alpha_i, \alpha_j \in \mathbb{R}$ such that $\va \alpha_i \vb + \va \alpha_j \vb >0$ and  $\Delta \pl \alpha_i \xi_{X_i} + \alpha_j \xi_{X_j} \pr =0$, where $\xi_A$ denotes the characteristic function on $A$. It is easy to check that
$\alpha_i \xi_{X_i} + \alpha_j \xi_{X_j}$ belongs to $T^n (C_{\mathbb{R}} (X))$ for every $n \in \mathbb{N}$, and consequently $T$ is not a shift.

In particular, we see that neither $C_{\mathbb{R}} \pl \mathbb{T}  + \mathbb{T}^{2} + \mathbb{T}^{3}   + \mathcal{N} \cup \{\pmb{\infty}\} \pr$  nor  $C_{\mathbb{R}} \pl \mathfrak{M} + \mathbb{T}  + \mathbb{T}^{2} + \mathbb{T}^{3}   + \mathcal{N} \cup \{\pmb{\infty}\} \pr$  admit an isometric shift.
\end{ex}

\subsection{Notation and assumptions}\label{pol}

\subsubsection{The numbers $p_n$, and the sets $\mathbb{P}_0$ and $\mathbb{N}_0$}\label{niko1}
We take $\mathbb{P}_0$ and $(p_n)_{n \in \mathbb{P}_0}$ to be an initial subset of $\mathbb{N}$ and a $\mathbb{P}_0$-compatible sequence, respectively. 
$\mathbb{N}_0$ will be taken as the set $\mathbb{N}$ if $\mathbb{P}_0 = \mathbb{N}$, and it will be taken as  $\{1,  \ldots, p_1   + \cdots + p_{P_0}\}$ if $\mathbb{P}_0 = \{1,  \ldots, P_0\}$.

\subsubsection{The sets $A_n$ and the maps $s_n$ and $\pi$}\label{niko2} Make sets on $\Bbb{N}_0$ in the following way. First $A_1 := \{1,  \ldots , p_1 \}$,  and in general, if $ n +1 \in \mathbb{P}_0$, $n \ge 1$, then take $A_{n+1} := \{p_1 + \cdots+ p_n   +1, \ldots , p_1 + \cdots+ p_n +p_{n+1} \}$.  Next, for each $n \in \mathbb{P}_0$, we put \[A_n = \{a_1^n , \ldots , a_{p_n}^n \}.\] 

Now define $s_n: A_n \ra A_n$ as $s_n (a_1^n) :=a_{p_n}^n$, and $s_n (a_j^n) :=a_{j-1}^n$, for $j =2, 3, \ldots, p_n$. 
Finally  define the map $\pi :\mathbb{N}_0 \ra \mathbb{P}_0$ sending each $k \in \mathbb{N}_0$ into the only number $n \in \mathbb{P}_0$ with $k \in A_n$.

\subsubsection{The basic family}\label{ajugar}
We assume that $ (Z_n, h_n, 1_n)_{n \in \mathbb{P}_0 }$ is $\mathbb{L}_0$-transitive, where $\mathbb{L}_0 := \{2 p_n : n \in \mathbb{P}_0\}$.  By the way we have taken the numbers $p_n$, it is obvious that when $\mathbb{P}_0:=\{1, \ldots, P_0\}$,     this fact is equivalent to saying that it is $\{2 p_{P_0}\}$-transitive.

\subsubsection{The spaces $X_n$ and $X_0$, and the maps $\chi_n$}\label{hn}

Define $X_n$ as the topological sum of $p_n$ spaces $Z_n$, that is, $X_n$ coincides (homeomorphically) with the  product $A_n \times Z_n$. 

Let us define next $\chi_n : A_n \times Z_n \ra A_n \times Z_n$ as $$\chi_n (a^n_j, z) := (s_n (a^n_j), h_n^{-1} (z)),$$ for every $j \in \{1,  \ldots, p_n\}$ and $z \in Z_n$.  

Our next step consists of defining the space $X_0$ as the topological sum of all $X_n$, which is clearly completely regular. A point $x \in X_0$ will be represented as $x=(k, z)$ where $k \in \mathbb{N}_0$ and $z \in Z_{\pi (k)}$.

\subsubsection{The maps $\phi$ and $a$ defined on $X_0$}\label{ht}
The complete definition of $\phi$ and $a$ will require more than one  step. The first one will be given here, and the others will be made later by extending them continuously. 
We start defining  $\phi : X_0 \ra X_0$ as $\phi (k, z) := \chi_{\pi (k)} (k, z)$ for every $(k,z) \in X_0$. It is easy to see that $\phi$ is a   homeomorphism. We also define the continuous function  $a: X_0  \ra \mathbb{K}$ as:
\begin{itemize}
\item $a (a_1^n, z) = -1$ for $n \in \mathbb{P}_0$ and $z \in Z_n$, 
\item $a \equiv 1$ on the rest of $X_0$.
\end{itemize}

\subsubsection{The spaces $X$ and $\omega X_0$,  and the complete definition of $\phi$ and $a$}\label{hk}
We consider $\omega X_0$ as any compactification of $X_0$ such that the following two properties hold:
\begin{itemize}
\item the map $\phi$ admits a continuous extension  to $\omega X_0$,  $\phi : \omega X_0 \ra \omega X_0$, which is a   homeomorphism.
\item the map $a$ admits a continuous extension  to  $\omega X_0$, $a : \omega X_0 \ra \mathbb{K}$.
\end{itemize}

It is clear that a space like this always exists, as we can take $\omega X_0$ equal to the Stone-\v{C}ech compactification of $X_0$. Nevertheless, for instance in Theorems~\ref{13dic22}, \ref{nesko},  and \ref{toallas}, as well as in Corollaries~\ref{kolya} and \ref{lanudo}, $\omega X_0$ will be a much simpler compactification. Notice that we are not assuming in principle that each $Z_n$ is compact, so a process of compactification is in general done even when $\mathbb{P}_0$ is finite.

Finally define our space $X$ as the topological sum of $\omega X_0$ and 
 $\mathcal{N}   \cup \{\pmb{\infty}\}$. Then extend  $a$ to a continuous map
 $a: X \setminus \{\mathbf{1}\} \ra \mathbb{K}$ by defining $a \equiv 1$ on
 $\mathcal{N}   \cup \{\pmb{\infty}\}$. This will be the complete definition of  $a$.

\subsubsection{The numbers $\gamma_n$}\label{hu}
For $n \in \mathbb{N}_0$, we consider any numbers  $\gamma_n \in \mathbb{K}$ satisfying the following properties: 

\begin{itemize}
\item There are at least two $n \in \mathbb{N}_0$ with $\gamma_n \neq 0$, 
\item $\sum_{n \in \mathbb{N}_0} \va \gamma_n \vb \le 1$, and
\item the determinant of the $p_n \times p_n$ matrix
\[M^{\gamma}_{n} := \left( 
\begin{array}{lllll}
\gamma_{a_1^n} & \gamma_{a_2^n} & \ldots & \gamma_{a_{p_n -1}^n } &\gamma_{a_{p_n}^n} \\
- \gamma_{a_{p_n}^n} & \gamma_{a_1^n} & \ldots &  \gamma_{a_{p_n -2}^n } &  \gamma_{a_{p_n -1}^n} \\
- \gamma_{a_{p_n -1}^n  } & - \gamma_{a_{p_n}^n} & \ldots &   \gamma_{a_{p_n -3}^n } &  \gamma_{a_{p_n -2}^n } \\
\vdots & \vdots &\ddots &\vdots &\vdots \\
- \gamma_{a_2^n} & - \gamma_{a_3^n} & \ldots & - \gamma_{a_{p_n}^n} &  \gamma_{a_1^n}
\end{array}
\right) \]
is different from $0$ for every $n \in \mathbb{P}_0$.
\end{itemize}

\subsubsection{The vectors $\mathbf{v}^n_i$, $\mathbf{f}_N^n$, and the numbers $\mathbf{v}^n_i \cdot \mathbf{f}_N^n$}\label{niko3}
Let $n \in \mathbb{P}_0$. For $i= 0, \ldots, p_n -1$, define the vector $\mathbf{v}_i^n$ as the $(i+1)$-th row of the matrix $M_n^{\gamma}$, that is, $$\mathbf{v}_0^n := \pl \gamma_{a^n_1}, \gamma_{a^n_2}, \ldots,\gamma_{a^n_{p_n -1}}, \gamma_{a^n_{p_n}} \pr,$$ and 
$$\mathbf{v}_i^n := \pl - \gamma_{a^n_{p_n -i+1}}, -\gamma_{a^n_{p_n - i+2}}, \ldots, -\gamma_{a^n_{p_n}}, \gamma_{a^n_{1}}, \ldots,  \gamma_{a^n_{p_n -i}} \pr $$ if $i=1, \ldots, p_n -1$. 
Also,  for $i=p_n , \ldots, 2 p_n -1$, define $\mathbf{v}_{i}^n := - \mathbf{v}_{i - p_n}^n$.

Given $f \in C (X)$ and $N \in \mathbb{N}$, we set $$\mathbf{f}_N^n  := \pl f \pl a_1^{n},  h_n^N \pl 1_n \pr \pr , \ldots , f \pl a_{p_n}^{n},  h_n^N \pl 1_n \pr \pr \pr.$$ 

Finally, $\mathbf{v}^n_i \cdot \mathbf{f}_N^n$ will denote the usual scalar product of $\mathbf{v}^n_i$ and  $\mathbf{f}_N^n$ in $\mathbb{K}^{p_n}$.

\subsubsection{The map  $\Delta$  and the isometry $T$}\label{conpermiso}

We now define $\Delta $ in the dual space $ C (X) '$ as
\[ \Delta (f)  := \sum_{n \in \mathbb{P}_0} \pl \sum_{m \in A_n} \gamma_m f (m, 1_n) \pr\]   for every $f \in C(X)$. Obviously $\vc \Delta \vd \le 1$.

Finally we introduce  a codimension $1$ linear isometry $T := T[a, \phi, \Delta]$. In particular we see that, for $k \in \Bbb{N}$ and $z \in Z_{\pi (k)}$, $(T f) (k, z) := - f( \phi(k, z))$ if $k$ coincides with  $a_1^{\pi(k)}$, and $(Tf) (k, z) =  f( \phi(k, z))$ otherwise. It is easy to check that $T$  is not an isometric  shift of type $\mathrm{II}$ since we deal with at least one $p_n \neq 1$ and two $\gamma_n \neq 0$.

\subsection{The results}

\begin{lem}\label{suso}
Assume that we follow the notation given in Subsection~\ref{pol}. If  $f \in \bigcap_{i=1}^{\infty} T^i (C(X))$, then $f$ is constant on $\mathcal{N} \cup \{\pmb{\infty}\}$. Moreover$$  f( \mathcal{N} \cup \{\pmb{\infty}\}) \equiv \sum_{n=1}^{\infty} \mathbf{v}^n_{[ k \mod 2 p_n]} \cdot \mathbf{f}_N^n ,$$ where this value is the same for every $N \in \mathbb{N}$ and $k \in \mathbb{N} \cup \{0\}$. 
\end{lem}

\begin{proof}
We will give the proof in the case when $\mathbb{P}_0 = \mathbb{N}$, since the other cases are similar.

It is easy to see that, if $k =t p_n$ for some $t \in \mathbb{Z}$ and $n \in \mathbb{N}$, then   $\phi^k (m, z) = (m,  h_n^{- k }(z))$ whenever $m \in A_n$ and $z \in Z_n$. From this we easily deduce the following claim.

\medskip

{\sc Claim 1.} {\em Let $m \in A_n$ and $k =t p_n$ for some $t \in \mathbb{Z}$, $n \in \mathbb{N}$. We have that, given any $z \in Z_n$,
\begin{itemize}
\item if $t$ is odd, then $(T^{k} f) (m, h_n^{k} (z)) = -f(m,  z)$, and
\item if $t$ is even, then $(T^{k} f) (m,  h_n^{k} (z)) = f(m, z)$.
\end{itemize}
}

\medskip

Now, the proof of the next claim is easy, taking into account that $f( \mathbf{k}) =   (T^{-k +1} f) (\mathbf{1})$.

\medskip

{\sc Claim 2.} {\em For every $\mathbf{k} \in {\mathcal N}$, 
\begin{eqnarray*}
f (\mathbf{k}) &=&  \Delta (T^{-k} f) \\&=&  \sum_{n \in \mathbb{P}_0} \pl
\sum_{m \in A_n} \gamma_m (T^{-k} f) (m, 1_n) \pr.
\end{eqnarray*}
}

\medskip

{\sc Claim 3.} {\em Let $n, N \in \mathbb{N}$, and  
 $k \in \mathbb{N} \cup \{0\}$. Then $$\sum_{m \in A_{n}} \gamma_m (T^{-k}  f )(m, h_{n}^{N-k}  (1_n)) =\mathbf{v}^n_{[k \mod 2p_n]} \cdot \mathbf{f}_N^n.$$ 
}

Let us prove the claim. For simplicity, we denote $y: = h_n^N (1_n)$ and    $w_i := f ( a_i^{n},  y )$, for $i=1, \ldots, p_n$, that is,  $\mathbf{f}_N^n  = \pl w_1, \ldots, w_{p_n} \pr$.
Notice first that, for every $(m, z) \in X_n$,  $f(m, z) = a(m, z) (T^{-1} f) (\chi_n (m,z))$, so taking into account that the maps $a$ and  $1/a$ coincide, we see that $(T^{-1} f) (m,z) = a(\chi_n^{-1} (m, z)) f(\chi_n^{-1} (m, z))$.  Consequently, if $m \in A_{n}$, then \[ (T^{-1} f)(m, h^{-1}_{n} (y)) = a (\chi_{n}^{-1} (m,  h^{-1}_{n}(y))) f (\chi_{n}^{-1} (m, h^{-1}_{n} (y))) ,\]
which means that  
$(T^{-1} f )(a_{p_n}^{n} , h_{n}^{-1} (y)) =  - f( a_{1}^{n},  y )$, and  $(T^{-1} f )(a_i^{n} , h_{n}^{-1} (y)) =  f ( a_{i+1}^{n},  y )$ for
$i \in \{1, 2, \ldots, p_{n} -1\}$. It is clear that in this case 
\begin{eqnarray*}
\sum_{i=1}^{p_{n}} \gamma_{a_i^{n}} (T^{-1}  f )(a_i^{n}, h_{n}^{-1}  (y))  &=&  - \gamma_{a_{p_n}^{n}} w_1  + \gamma_{a_1^{n}} w_2   + \ldots +   \gamma_{a_{p_{n} -1}^{n}} w_{p_{n}}\\
&=& \mathbf{v}^n_{1} \cdot \mathbf{f}_N^n.
\end{eqnarray*}

Of course a similar process shows that, if $k \in \{2, \ldots, p_n -1\}$, then
\begin{eqnarray*}
(T^{-k} f )(m, h_{n}^{-k} (y)) &=& a (\chi_{n}^{-1} (m,  h_{n}^{-k} (y))) (T^{1-k} f) (\chi_{n}^{-1} (m, h_{n}^{-k} (y) )) \\
&=&  \pl \prod_{i=1}^2 a (\chi_{n}^{-i} (m,  h_{n}^{-k} (y)))  \pr  (T^{2-k} f) (\chi^{-2}_{n} (m,  h_{n}^{-k} (y)))\\
&\vdots&
\\
&=& \pl \prod_{i=1}^k a (\chi_{n}^{-i} (m,  h_{n}^{-k} (y)))  \pr    f (\chi_{n}^{-k} (m , h_{n}^{-k} (y))),
\end{eqnarray*}
that is, 
\begin{eqnarray*}
(T^{-k} f )(a_{p_n}^{n} , h_{n}^{-k} (y)) &=&  - w_{k},\\
(T^{-k} f )(a_{p_n -1}^{n} ,h_{n}^{-k} (y)) &=&  - w_{k -1}, \\
&\vdots& \\
(T^{-k} f )(a_{p_n -k +1}^{n} , h_{n}^{-k} (y)) &=&  - w_{1}, \mbox{ and} \\
 (T^{-k} f )(a_j^{n} ,h_{n}^{-k} (y)) &=&  w_{k+j}
\end{eqnarray*}
 for
$j \in \{1,  \ldots, p_{n} -k\}$. As above this implies that
\begin{equation}\label{lila}
\sum_{m \in A_n} \gamma_m (T^{-k} f )(m,h_{n}^{-k} (y)) = \mathbf{v}^n_{k} \cdot \mathbf{f}_N^n.
\end{equation}

Next, suppose for instance that $k= ap_n +b$, where $a$ is odd and $b \in \{0, \ldots , p_n -1\}$. According to Claim 1, 
we have that
\begin{eqnarray*}
(T^{-ap_{n} -b } f) (m, h_{n}^{-ap_{n} - b} (y)) &=& (T^{-ap_{n}} (T^{-b} f)) (m, h_{n}^{-ap_{n} - b} (y)) \\
&=& - (T^{-b} f )(m, h_{n}^{-b} (y)),
\end{eqnarray*}
and it is immediate from Equality~\ref{lila} that 
\[\sum_{m \in A_n} \gamma_m (T^{-ap_{n} - b}  f )(m, h_{n}^{-ap_{n} - b}  (y)) = - \mathbf{v}^n_{b} \cdot \mathbf{f}_N^n = \mathbf{v}^n_{p_n + b} \cdot \mathbf{f}_N^n.\]

The case when $a$ is even is similar. This finishes the proof of the claim.

\medskip

Next fix any $\epsilon >0$, and
take $n_0 \in \mathbb{N}$  such that $$\vc f \vd \sum_{n= n_0 +1}^{\infty} \pl \sum_{m \in A_n} \va \gamma_m \vb \pr  < \epsilon.$$
Define, for each $k \in \mathbb{N}$,
$$B (k) := \sum_{n=1}^{n_0} \pl \sum_{m \in A_n} \gamma_m (T^{-k} f) (m,  1_n) \pr.$$
 
By Claim 2, it is clear that for every $k$,
\begin{equation}\label{tarlow}
\va f(\mathbf{k}) -  B(k) \vb < \epsilon.
\end{equation}

Now  make a partition of the set of natural numbers as follows. For $i = 0, \ldots,  2p_{n_0} -1$, let $$\mathbb{P}_i :=\{n \in \mathbb{N}: n = i \mod 2 p_{n_0}\}.$$  

Fix any $N \in \mathbb{N}$. Since $(Z_n, h_n, 1_n)_{n \in 
\mathbb{P}_0}$  is  $\mathbb{L}_0$-transitive and  $2 p_{n_0} \in \mathbb{L}_0$, then  $(h_1^N (1_1), \ldots, 
h_{n_0}^N (1_{n_0})) $ is a limit point  of $\{(h_1^l (1_1), \ldots, h_1^l (1_{n_0})) : l \in \mathbb{P}_i\}$ in 
$ \prod_{n=1}^{n_0} Z_n$, for each $i=0, \ldots, 2 p_{n_0} -1$. Then, taking into account that the map
 $P_f: \prod_{n=1}^{n_0} Z_n  \ra \mathbb{K}^{p_1 + \cdots + p_{n_0}}$ 
sending each $(z_1, \ldots, z_{n_0})$ into 
$$P_f (z_1, \ldots, z_{n_0}) := (f(a_1^1, z_1) , \ldots, f(a_{p_1}^1, z_1), \ldots, f(a_1^{n_0}, z_{n_0}) , \ldots,
 f(a_{p_{n_0}}^{n_0}, z_{n_0}))$$ is continuous, we deduce that there exists an increasing subsequence $(u_j^i )$ in $\mathbb{P}_i$
 such that, for  
 each $n \in \{1, 2, \ldots, n_0\}$ and  each $m \in A_n$, 
\[\lim_{j \ra \infty} f(m , h_n^{u_j^i} (1_n)) = f(m, h_n^N (1_n)) .\]

\medskip
 
   Consider each of the sequences $(u_j^i)$ given above.  By Claim 3, it is clear that due to our assumptions on the numbers $p_n$,
\begin{eqnarray*}
B (u_j^i) &=& \sum_{n=1}^{n_0} \mathbf{v}^n_{[  u_j^i \mod 2 p_n]} \cdot \mathbf{f}^n_{u^i_j} \\
&=& \sum_{n=1}^{n_0} \mathbf{v}^n_{[  i \mod 2 p_n]} \cdot \mathbf{f}^n_{u^i_j}. 
\end{eqnarray*} 
Now,   using Equality~\ref{tarlow} and taking limits   we obtain that, for each $i=0,  \ldots  , 2p_{n_0} -1$, 
$$  \va f( \pmb{\infty}) - \sum_{n=1}^{n_0} \mathbf{v}^n_{[ i \mod 2 p_n]} \cdot  \mathbf{f}_N^n \vb \le \epsilon.$$
We  conclude that, for every $k \in \mathbb{N} \cup \{0\}$ and $N \in \mathbb{N}$,
$$  f( \pmb{\infty}) = \sum_{n=1}^{\infty} \mathbf{v}^n_{[ k \mod 2 p_n]} \cdot  \mathbf{f}_N^n.$$

In particular, the sum $\sum_{n=1}^{\infty} \mathbf{v}^n_{[ k \mod 2 p_n]} \cdot  \mathbf{f}_N^n$ does not depend neither on $k$ nor on $N$. On the other hand, it is clear by Claims 2 and 3 that, taking $N=k$,  $$  f( \mathbf{k}) = \sum_{n=1}^{\infty} \mathbf{v}^n_{[ k \mod 2 p_n]} \cdot \mathbf{f}_N^n $$ for every $\mathbf{k} \in \mathcal{N}$, so we deduce that $f$ is constant on $\mathcal{N} \cup \{\pmb{\infty}\}$.
\end{proof}

\begin{lem}\label{gotera}
Assume that we follow the notation given in Subsection~\ref{pol}. Then the  set $\tl \chi_n^k (a_1^{n}, 1_n) : k \in \mathbb{Z} \tr$ is dense in  $X_n$ for each $n \in \mathbb{P}_0$.
\end{lem}

\begin{proof}
First, for $n \in \mathbb{P}_0$ fixed, $(Z_n, h_n, 1_n)$ is transitive by hypothesis, so $\tl \chi_{n}^{-k} \pl a^n_1,  h_n^i (1_n) \pr : k \in p_n \mathbb{N}, i=0, \ldots, p_n-1 \tr$  is dense in $\{a_1^n\} \times Z_n$.
Also, since $(Z_m, h_m, 1_m)_{m \in \mathbb{P}_0}$ is $\{2p_n\}$-transitive, it is easy to check that, for 
$i =1, \ldots, p_n -1$, each point $h_n^i (1_n)$ belongs to the closure of $\tl h_n^{2jp_n } (1_n) :j \in \mathbb{N}\tr$. Taking  both facts into account, we see that $\tl \chi_{n}^{-k} \pl a^n_1,  1_n \pr : k \in p_n \mathbb{N}\tr$  is dense in $\{a_1^n\} \times Z_n$, and  this implies that $\tl \chi_{n}^{-k} \pl a^n_1,  1_n \pr : k \in  \mathbb{N}\tr$  is dense in $X_n$.
\end{proof}

\begin{thm}\label{susovict}
Assume that we follow the notation given  in Subsection~\ref{pol}.  Then $T$ is a shift  of type $\mathrm{I}_0$. Also, $T$ is $\infty$-generated if $\mathbb{P}_0 = \mathbb{N}$, and it is $P_0$-generated if $\mathbb{P}_0 = \{1, 2, \ldots, P_0\}$.
\end{thm}

\begin{proof}

We split the proof into two parts.

\medskip

{\em Case 1. Suppose  that $\mathbb{P}_0 =\{1\}$.} Recall that we are assuming in this case that $p_1 \neq 1$. 
Let us prove that, if $f \in \bigcap_{l=1}^{\infty} T^l (C(X))$, then  $f=0$. Fix any $i =0, \ldots, p_1 -1$. First we have, by Lemma~\ref{suso}  that, 
given any $N \in \mathbb{N}$,  $\mathbf{v}^1_{[ i \mod 2 p_1]} \cdot \mathbf{f}_N^1 = \mathbf{v}^1_{[ i +  p_1 \mod 2 p_1]} \cdot \mathbf{f}_N^1 $. Taking into account the definition of $\mathbf{v}^1_i$,
this means that   $- \mathbf{v}^1_{i} \cdot \mathbf{f}_N^1 = \mathbf{v}^1_{i} \cdot \mathbf{f}_N^1 $, and consequently the product of the matrix $M_1^{\gamma}$ and  the vector $\mathbf{f}_N^1$ is $\mathbf{0}$. Since  $M_1^{\gamma}$ is invertible, we conclude that $\mathbf{f}_N^1 = \mathbf{0}$. Since this happens for every $N \in \mathbb{N}$, we conclude by denseness that $f \equiv 0$ on $\omega X_0$. Finally the value given in Lemma~\ref{suso} for $f$ on $\mathcal{N} \cup \{\pmb{\infty}\}$ must also be $0$, and we are done.

The fact that $T$ is $1$-generated follows from Lemma~\ref{gotera}. 

\medskip

{\em Case 2. Suppose  that $\mathbb{P}_0 \neq \{1\}$.} 
We just give the proof in the case $\mathbb{P}_0 = \mathbb{N}$, since the other cases are similar.  Let $f \in \bigcap_{l=1}^{\infty} T^l (C(X))$ with $\vc f \vd \le 1$ . 

Let us prove first that $f \equiv 0$ on $\mathcal{N} \cup \{\pmb{\infty}\}$. 
By hypothesis we have that if $n <i$, then $p_i = 2 l p_n$ for some $l \in \mathbb{N}$, so  
taking into account that $\mathbf{v}^n_{[ k \mod 2 p_n]} = - \mathbf{v}^n_{[p_n + k \mod 2 p_n]}$ for every $k \in \mathbb{N} \cup \{0\}$,
we deduce that       
\begin{eqnarray*}
\sum_{k=1}^{p_i} \mathbf{v}^n_{[ k \mod 2 p_n]}  &=& \sum_{k=1}^{2 p_n} \mathbf{v}^n_{[ k \mod 2 p_n]}  + \cdots + \sum_{k= 2 (l-1) p_n +1}^{2 l p_n} \mathbf{v}^n_{[ k \mod 2 p_n]} \\
&=& 0 .
\end{eqnarray*}

This implies that,  for every $N \in \mathbb{N}$ and  $i \in \mathbb{N} \setminus \{1\}$,
\begin{eqnarray*}
\sum_{k=1}^{p_i} f(\mathbf{k}) &=& \sum_{k=1}^{p_i} \pl \sum_{n=1}^{\infty} \mathbf{v}^n_{[ k \mod 2 p_n]} \cdot \mathbf{f}_N^n \pr \\
&=& \sum_{k=1}^{p_i} \pl \sum_{n=1}^{i-1} \mathbf{v}^n_{[ k \mod 2 p_n]} \cdot \mathbf{f}_N^n   + \sum_{n=i}^{\infty} \mathbf{v}^n_{[ k \mod 2 p_n]} \cdot \mathbf{f}_N^n  \pr \\
&=&  \sum_{n=1}^{i-1} \sum_{k=1}^{p_i}  \mathbf{v}^n_{[ k \mod 2 p_n]} \cdot \mathbf{f}_N^n   + \sum_{k=1}^{p_i} \sum_{n=i}^{\infty} \mathbf{v}^n_{[ k \mod 2 p_n]} \cdot \mathbf{f}_N^n  \\
&=& \sum_{k=1}^{p_i} \sum_{n=i}^{\infty} \mathbf{v}^n_{[ k \mod 2 p_n]} \cdot \mathbf{f}_N^n .
\end{eqnarray*}

Now, by Claim 3 in Lemma~\ref{suso}, and taking into account that $f$ is constant on $\mathcal{N} \cup \{\pmb{\infty}\}$, we have that $$p_i \va f(\mathbf{1}) \vb \le \sum_{k=1}^{p_i} \sum_{n=i}^{\infty} \va \mathbf{v}^n_{[ k \mod 2 p_n]} \cdot \mathbf{f}_N^n  \vb \le p_i \sum_{n=p_1 + \cdots + p_{i-1} + 1}^{\infty} \va \gamma_n \vb ,$$ for every $i \in \mathbb{N} \setminus \{1\}$. Consequently, as $i$ goes to infinity, the term $ \sum_{n=p_1 + \cdots + p_{i-1} + 1}^{\infty} \va \gamma_n \vb$ goes to $0$, and we easily conclude that $f(\mathbf{1})=0$, and  $f \equiv 0$  on $\mathcal{N} \cup \{\pmb{\infty}\}$.

\medskip

We deduce then that $\sum_{n=1}^{\infty} \mathbf{v}^n_{[ k \mod 2 p_n]} \cdot \mathbf{f}_N^n =0$ for every $k \in \mathbb{N} \cup \{0\}$ and $N \in \mathbb{N}$.

We must now show that this implies that $\mathbf{f}_N^n = \mathbf{0} \in \mathbb{K}^{p_n}$. First notice that for $n_0  \in \mathbb{N} \setminus \{1\}$ and $k \in \mathbb{N} \cup \{0\}$, 

\begin{eqnarray}\label{savor}
 \va 2  \mathbf{v}^{n_0}_{[k \mod 2 p_{n_0}]} \cdot \mathbf{f}_N^{n_0} \vb &=& \va \mathbf{v}^{n_0}_{[k \mod 2 p_{n_0}]} \cdot \mathbf{f}_N^{n_0} -  \mathbf{v}^{n_0}_{[ k + p_{n_0} \mod 2 p_{n_0}]} \cdot \mathbf{f}_N^{n_0} \vb \nonumber \\ &=& \va \sum_{n=1}^{n_0} \mathbf{v}^n_{[k \mod 2 p_n]} \cdot \mathbf{f}_N^{n} - \sum_{n=1}^{n_0} \mathbf{v}^n_{[k+ p_{n_0} \mod 2 p_n]} \cdot \mathbf{f}_N^{n}\vb ,
 \end{eqnarray}
 and 
 \begin{eqnarray}\label{sabor}
 \va 2 \sum_{n=1}^{n_0 -1} \mathbf{v}^n_{[k \mod 2 p_n]} \cdot \mathbf{f}_N^{n}  \vb 
 &=& \va \sum_{n=1}^{n_0 -1} \mathbf{v}^n_{[k \mod 2 p_n]} \cdot \mathbf{f}_N^{n} + \sum_{n=1}^{n_0 -1} \mathbf{v}^n_{[k+  p_{n_0} \mod 2 p_n]} \cdot \mathbf{f}_N^{n} \vb \nonumber \\ 
&=& \va \sum_{n=1}^{n_0} \mathbf{v}^n_{[k \mod 2 p_n]} \cdot \mathbf{f}_N^{n}+ \sum_{n=1}^{n_0} \mathbf{v}^n_{[ k+ p_{n_0} \mod 2 p_n]} \cdot \mathbf{f}_N^{n} \vb .
 \end{eqnarray}

Next take  any $\epsilon >0$ and $n_1 \in \mathbb{N}$ 
such that $\sum_{n=n_1+1}^{\infty}  \pl \sum_{m \in A_n} \va \gamma_m  \vb\pr < \epsilon$. 
This obviously gives that $\va \sum_{n=1}^{n_1} \mathbf{v}^n_{[ k \mod 2 p_n]} \cdot \mathbf{f}_N^n \vb < \epsilon$ for every $k \in \mathbb{N} \cup \{0\}$ and $N \in \mathbb{N}$. Let us prove that $\va \mathbf{v}^n_{[ k \mod 2 p_n]} \cdot \mathbf{f}_N^n \vb < \epsilon$ for every $n \le n_1$. Of course the result is true for the case when $n_1 =1$. Next, if we  assume that it holds for $n_1 \le l$, we  can easily prove it for $n_1 =l+1$ using Equalities~\ref{savor} and \ref{sabor}.
 
We deduce that $\mathbf{v}^n_{[ k \mod 2 p_n]} \cdot \mathbf{f}_N^n=0$ for every $n, N \in \mathbb{N}$ and 
$k \in \mathbb{N} \cup \{0\}$. Due to the properties of each  matrix $M_n^{\gamma}$, we conclude that $\mathbf{f}_N^n = \mathbf{0}$ and, by denseness, $f$ must be equal to $0$, and $T$ is a shift. 

Finally, notice that, given a point $x \in X_n$, the set $\{\phi^{k} (x) : k \in \mathbb{Z}\}$ is contained in $X_n$, which implies that $T$ is $\infty$-generated. Notice also that if $\mathbb{P}_0 \neq \mathbb{N}$, to prove that $T$ is $P_0$-generated we use both the above  fact and Lemma~\ref{gotera}.
\end{proof}

\begin{rem}
In the proof of the previous theorem, we borrow some of the ideas of \cite[Theorem 3.5]{GN}, where the sequence $(2^n)$  plays a fundamental r\^{o}le when proving the existence of an isometric shift of type $\mathrm{I}$ on $C(X)$, for $X= \beta \mathbb{N} + \mathcal{N} \cup 
\{\pmb{\infty}\}$. On the other hand, we must add that  the proof of Lemma~\ref{suso} is valid even if  we just assume that the sequence $(p_n)_{n \in \mathbb{P}_0}$ we take satisfies that each $p_{n+1}$ is an integer  multiple of $p_n$
 (and  $p_{n+1} > p_n$), instead of being an  $\mathbb{P}_0$-compatible sequence.  Nevertheless, in Theorem~\ref{susovict} our requirements on the sequence $(p_n)$ are necessary. Let us see it:  Suppose that there exist two different $p_n$ (for instance, $p_1$ and $p_2$) which are odd, and  that
  $\gamma_1 = \gamma_{p_1 +1} \neq 0$, and  $\gamma_i =0$ for every other $i \le p_1 + p_2$. Consider the function $f \in C (X)$ satisfying
\begin{itemize}
\item $f(\{a_i^1\} \times Z_1) \equiv 1$ if $i \in \{1, \ldots, p_{n_1}\}$ is odd, and $f(\{a_i^1\} \times Z_1) \equiv - 1$ 
if $i$ is even, and
\item $f(\{a_i^2\} \times Z_2) \equiv - 1$ if $i \in \{1, \ldots, p_{n_2}\}$ is odd,
 and $f(\{a_i^2\} \times Z_2) \equiv  1$ if $i$ is even,
 \end{itemize}
and is equal to $0$ in the rest of $X$. It is easy to check that $\sum_{n=1}^{\infty} \mathbf{v}^n_{[ k \mod 2 p_n]} \cdot \mathbf{f}_N^n \equiv 0$ for every $k \in \mathbb{N} \cup \{0\}$ and every $N \in \mathbb{N}$. Nevertheless we can see that $Tf = -f$, and consequently $T$ is not a shift.
\end{rem}

\begin{rem}\label{notransi}
Let us see that the condition  on  
$ (Z_n, h_n, 1_n)_{n \in \mathbb{P}_0 }$ of 
being $\mathbb{L}_0$-transitive (following the notation given in Subsection~\ref{ajugar}) is not redundant in  Theorem~\ref{susovict}. Consider for instance the case of $Z_1$ 
consisting of the union of 
two copies of a totally transitive set $W$ (with respect to a homeomorphism $\mathbf{r}$ and a point $w$), that is, 
$Z_1:=  \pl \mathbb{Z}/2 \mathbb{Z} \pr \times W $. It is  clear that the point $1_1 := (\overline{0},w) \in \pl \mathbb{Z}/2 \mathbb{Z} \pr \times W$ is transitive for the map $h_1 : Z_1 \ra Z_1$, defined as 
$h_1 (\overline{a}, x) := (\ \overline{a+1}, \mathbf{r} (x))$  
for every  $(\overline{a}, x) \in \pl \mathbb{Z}/2 \mathbb{Z} \pr \times W$.
Let $A_1 := \{1, 2\}$ and $X_1 := \{1,2\} \times Z_1$. We may also consider some other sets $X_n$ defined as in Subsection~\ref{pol}, although they will play no r\^{o}le in what follows. Define finally $\chi_1: X_1  \ra X_1$, $X$, $\phi$ and $T$ using the  notation given in Subsections~\ref{hk} and \ref{conpermiso}. In particular, for $\mathbf{1} \in \mathcal{N}$,  $$(Tf) (\mathbf{1}) = \gamma_1 f(1, 1_1) + \gamma_2 f(2, 1_1)$$ for every $f \in C(X)$. We have that $T$ is a codimension $1$ 
linear isometry, but it is not a shift because the function $f$ defined as $f \equiv \gamma_1$ on $ \{2\} \times \{\overline{0}\} \times W  $, $f \equiv - \gamma_2$ on $ \{1\} \times \{\overline{0}\} \times W  $, and $f \equiv 0$ on the rest of $X$,  satisfies $T^2 f = - f$.
\end{rem}
 
\subsection{Some consequences and  proofs}

In Subsection~\ref{hk}, we said that the compactification $\omega X_0$ could always be taken as the Stone-\v{C}ech compactification of $X_0$. In the proofs of Theorems~\ref{13dic22} and \ref{nesko} we will take advantage of the fact that we  deal with copies of compact spaces, and use appropriate one-point compactifications of some sequences of
them.

\begin{proof}[Proof of Theorem~\ref{13dic22}]
Consider the  complex linear space $\mathbf{B} := \ell^1 (\mathbb{R})$ of all $\mathbb{C}$-valued maps $f$ on $\mathbb{R}$  such that $f(x) \neq 0$ for at most countably many $x \in \mathbb{R}$, and such that  $\vc f \vd_1 := \sum_{x \in \mathbb{R}} \va f(x) \vb < + \infty$. $\vc \cdot \vd_1$ determines a Banach space structure on $\mathbf{B}$, whose dual $\mathbf{B}'$ can be identified with the 
(complex) Banach  space $ \ell^{\infty} (\mathbb{R})$ of all $\mathbb{C}$-valued bounded maps on $\mathbb{R}$ endowed with the supremum norm (see \cite[p. 137]{K}).

We have that   for every  $x \in \mathbb{R}$ and $b= 0, 1$, the set $A_x^b := \{ f \in \mathbf{B}' : \va f(x) \vb = b \}$ is   closed  in the $\mathrm{weak}^*$-topology. Consequently, if  $\mathbb{R}_n \subset \mathbb{R}$ is a set of cardinality $\kappa_n$ for each
  $n \in \mathbb{M}$, then 
 $C_n := \{f \in \mathbf{B}' : \va f(x) \vb =1 \hspace{.03in} \forall x \in \mathbb{R}_n, f(x) =0 \hspace{.03in}  \forall x \notin \mathbb{R}_n \}$ is a closed subset of $\mathbf{B}'$ in the $\mathrm{weak}^*$-topology; also, it is norm-bounded 
  in $\mathbf{B}'$, so by Alaoglu's Theorem it is compact. 
On the other hand, since each $C_n$ is homeomorphic to $\mathbb{T}^{\kappa_n}$, then  we can embed homeomorphically
 each $\mathbb{T}^{\kappa_n}$ in any open ball of $\mathbf{B}'$ (with the $\mathrm{weak}^*$-topology).
  
Fix a line $L_0$ passing through the origin $\mathbf{0} \in \ell^{\infty} (\mathbb{R})$, and a  sequence $(y_k)_{k \in \mathbb{N}}$ 
of pairwise distinct points in $L_0$ convergent to $\mathbf{0}$, in such a way that all of them lie on the same semiline
 determined by the origin, and that $\vc y_k \vd =1/k $ for all $k \in \mathbb{N}$. Consider also the sequence $(-y_k)_{k \in \mathbb{N}}$. 
Next, for 
each $\pm y_k$, take a line $L_{\pm k}$ passing through it, and being perpendicular to $L_0$ (for instance, assuming that all $L_{\pm k}$ and $L_0$ are located in the same plane). Obviously, each point $\pm y_k$ 
divides the line $L_{\pm k}$ into two semilines. We select one of them, and call it $S_{\pm k}$.

For each $\pm k$, we will consider a different sequence, starting at an index $I_k$ defined as follows: If $k= \pm 1, \ldots, \pm p_2/2$, then $I_k:= 2$, and in general, if $n \ge 2$ and $k = \pm (p_{n}/2 +1), \ldots, \pm p_{n+1}/2$, then $I_k := n+1$. Now for each $\pm k$, consider a sequence of pairwise disjoint open balls, which will be indexed
$\pl B^{\pm k}_n \pr_{n \ge I_k}$, with radii decreasing to $0$, and with centers in 
$S_{\pm k}$ 
converging to $y
_{\pm k}$. 
We can do it so   that  all balls in the set $\tl B^{\pm k}_n : k \in \mathbb{N},  n \ge I_k \tr$ are pairwise disjoint, and each ball  $B^{k}_n$ intersects just the line $L_k$ and no other (not even $L_0$), for every   $k \in \mathbb{N} \cup - \mathbb{N}$ and $ n \ge I_k$. Notice that according to the way we have indexed the balls, we have that, for $n \ge 2$, the only balls having $n$ as a subindex  are $B_n^k$ for $k=\pm 1  , \ldots, \pm p_n/2$. Define finally some $p_1$ new open  balls having $1$ as a subindex, namely $B_1^k$  for $k=1, \ldots, p_1$, which we suppose are disjoint from all the other balls taken above.
Now, for each $n \in \mathbb{P}$, consider $i \in \mathbb{M}$ such that $p_n$ belongs to $\mathscr{A}_i$, and 
call $Z_n := \mathbb{T}^{\kappa_i}$. Next, by Example~\ref{zarrau}, for each $n \in \mathbb{P}$, we can find a homeomorphism $h_n : Z_n \ra Z_n$, and a point $1_n \in  Z_n$ such that $\pl Z_n, h_n, 1_n \pr$ is $\mathbb{N}$-transitive. Obviously the case $\kappa_i= 0$, that is, when the space we deal with is the single point $\mathbb{T}^0$,  does not pose any problems.

Next put a copy of $Z_n$ inside each ball $B^{k}_n $ (for $k= 1, \ldots, p_1$ if $n =1$, and $k = \pm 1  , \ldots, \pm p_n/2$ otherwise), and denote this 
corresponding copy by  $Y_n^k$. We see that we could be dealing just with a finite number of balls, if $\mathbb{M}$ is finite.

We are going to use the notation given in  Subsection~\ref{pol}, and then apply   Theorem~\ref{susovict}. 
We consider $\{a_1^1\} \times Z_1, \ldots, \{a_{p_1}^1\} \times Z_1$, and identify them with $Y_1^1, \ldots, Y_1^{p_1}$,
respectively. On the other hand, for each $n \ge 2$, $n \in \mathbb{P}$, we consider $\{a_1^n\} \times Z_n, \ldots, \{a_{p_n}^n\} \times Z_n$, and identify them with $Y_n^{-p_n/2}, \ldots, Y_n^{p_n/2}$, respectively. We see that
the set $Y_0$ formed by the union of all $Y_n^k$, and the set $X_0$ formed by the topological sum of all $\{a_j^n\} \times Z_n$
(following the notation of Subsection~\ref{pol}) are homeomorphic. 

Consequently, if $\mathbb{P}$ is finite, Theorem~\ref{susovict} provides us a way to construct an isometric shift
as we want.

Let us finally study the case when $\mathbb{P}$ is not finite. Let us see that 
$Y_C:= Y_0 \cup \{\mathbf{0}\} \cup 
\{ y_{\pm k} : k \in \mathbb{N}\}$ is a compactification of $Y_0$. Obviously, we just need to prove that
$Y_C$ is the closure of $Y_0$ in $\mathbf{B}'$ with the $\mathrm{weak}^*$-topology. Notice that, for every $\epsilon>0$, we can cover the set $M:= \{\mathbf{0}\} \cup 
\{ y_{\pm k} : k \in \mathbb{N}\}$ with a finite number of closed balls $B_1, \ldots, B_m$ of radius $\epsilon$ centered at points of $M$.
Obviously all sets $Y_n^k$ but a finite number of them are contained in $B_1 \cup \ldots B_m$, which again by Alaoglu's Theorem is a compact set, and consequently closed. This implies that
the closure of $Y_0$ in the $\mathrm{weak}^*$-topology is contained in $Y_0 \cup B_1 \cup \ldots B_m$. Since this works for every $\epsilon >0$, we easily conclude that $Y_C$ is the closure of $Y_0$. This can be translated into a closure of $X_0$ by a countable family of points. It is clear that the conditions given in Subsection~\ref{pol} are satisfied, and Theorem~\ref{susovict} can be applied.
\medskip

As for the case when  $\kappa_n \le \aleph_0$ for all $n \in \mathbb{M}$, we have two possibilities. First it is clear that if 
$\mathrm{s} \in \mathbb{N}$, then we can do the above process in $\mathbb{C}^\mathrm{s}$. On the other hand, if we suppose that 
$\mathrm{s} = \aleph_0$, we can use the famous result that states that $\ell^2$ is homeomorphic to $\mathbb{R}^{\mathbb{N}}$ (see \cite{An}), which is
to say to $\mathbb{C}^{\mathbb{N}}$. As a consequence, a copy of $\mathbb{T}^{\kappa_n}$ can be embedded in each open ball
of $\ell^2$ (endowed with the norm topology), and we  do again  the same process as above, taking into account that $Y_C$ is now a compactification of $Y_0$ with respect to the norm topology.
\end{proof}

We  easily deduce the following corollary, using the case when $\mathrm{s}=1$.

\begin{cor}\label{kolya}
 In  $\mathbb{R}^2$, we can find a  compact set $X$ having a countably infinite number of components (each of them being infinite), and such that   $C(X)$ admits an isometric shift which is $\infty$-generated. 
\end{cor}

\begin{cor}\label{lanudo}
Let 
$E$ be a normed space, and let 
$A \subset E$ be compact and totally transitive (with respect to a homeomorphism). For each $n \in \mathbb{N} \cup \{\infty\}$, there exist a space  $X_A \subset E$ consisting of a topological sum of  copies of $A$, and a compactification $X_C$ of $X_A$ in $E$ such that, if  $(x_k)$ is an infinite sequence in $E \setminus X_C$ converging to a point $x_0 \notin X_C$, then for $X:= X_C \cup \{x_k : k \in \mathbb{N} \} \cup \{x_0\}$, $C(X)$ admits an isometric shift  of type $\mathrm{I}_0$ which is $n$-generated. Moreover 
$X_C$ does not have  isolated points if $A$ has at least two points.
\end{cor}

\begin{proof}
The proof of the first part is similar
 to that of Theorem~\ref{13dic22}, where we identify a sequence and its limit with $\mathcal{N} \cup \{\pmb{\infty}\} $. Also, if $n \in \mathbb{N}$, then $X_A$ consists of a finite number of copies of $A$, and coincides with $X_C$. Notice that if we are in the particular case that $E= \mathbb{K} = \mathbb{R}$, then the process is slightly different, since when we deal with infinitely many copies of  $A$ converging to $\pm y_k$, we must include them in the segment between $\pm y_{k+1}$ and $\pm y_k$.  On the other hand, we just mention that  when $A$ has at least two points, then the fact that it is totally transitive implies that it has no isolated points. Now it is easy to see that in that case $X_C$ has no isolated points either.
\end{proof}

\begin{rem}\label{manyanadel}
In Corollary~\ref{lanudo}, there are no restrictions on the dimension of $E$ (as far as it is nontrivial). Also the topology that we assume in $X_C$ is the inherited from the metric in  $E$. 
\end{rem}

The idea of the proof of Theorem~\ref{13dic22} and Corollary~\ref{lanudo} can be also used to prove Theorem~\ref{nesko}.

\begin{proof}[Proof of Theorem~\ref{nesko}]
Let us prove the first part. Notice that  we cannot follow exactly the same pattern as in the proof of Theorem~\ref{13dic22}, since it is not possible in general to embed each $K_n^{\kappa_n}$ in the Banach space $\mathbf{B}'$.  
Anyway we may consider  some spaces $Y^k$ (defined below), and work as in the previous proof. We do the process just 
for $\mathbb{P} = \mathbb{N}$ since other cases are similar.

For each $k \in \mathbb{Z} \setminus\{0\}$, we take $I_k$ as in the proof of Theorem~\ref{13dic22}, and put $Y^k_n := K_i^{\kappa_i}$ for each $n \ge I_k$, $n \in \mathscr{A}_i$ (where $\mathscr{A}_i$ is given as in Theorem~\ref{13dic22}). Next  define $Y^k$ as the topological sum of all $Y^k_n$, which is locally compact. Then  consider its compactification by one point, denoted $\infty_k$. 

Let $Y_0 := \bigcup_{\va k \vb=1}^{\infty} Y^k$. It is clear that if we identify $\infty_k$ with the number 
$1/k \in \mathbb{R}$ for each $k \in \mathbb{Z} \setminus \{0\}$, and in $ \{0\} \cup 
\{ \pm 1/k : k \in \mathbb{N}\}$ we consider the  usual topology induced from $\mathbb{R}$, then  $Y_C:= Y_0 \cup \{0\} \cup 
\{ \pm 1/k : k \in \mathbb{N}\}$ becomes a compactification of $Y_0$. It is also clear that a similar process as in the proof of Theorem~\ref{13dic22} relates this compactification  to the notation of Subsection~\ref{pol}. Then we finish by applying Theorem~\ref{susovict} (in a similar way as in Theorem~\ref{13dic22}) through Remark~\ref{zien}.

\medskip

Finally when $\kappa_n = \aleph_0$ and $K_n$ is metrizable for all $n$, then so is $K_n^{\mathbb{N}}$, which can  consequently  be embedded in 
$[0,1]^{\mathbb{N}}$ (see \cite[Theorem 23.1]{W}). As in the previous proof, each $K_n^{\mathbb{N}}$ can then be embedded in the unit ball of $\ell^2$, and the conclusion can be reached in a similar way as there.
\end{proof}

\begin{rem}\label{cincox}
In the first part of the proof of Theorem~\ref{nesko}, we can take any compactification $\omega Y_0$ of $Y_0$, as long as it satisfies the two requirements given in Subsection~\ref{hk}, and the proof remains valid.

\end{rem}

\section{A  method for obtaining examples of composition operators}\label{luzia}

In this section, we will analyze in particular the case when the weight $a \in C (X \setminus \{\mathbf{1}\})$  is constantly equal to $1$. A first approach shows that constructions similar to that given above (that is, based on the topological sum of a compact space and ${\mathcal N} \cup \{\pmb{\infty}\}$) do not yield an isometric shift.

\begin{prop}\label{nde}
Let $X$ be compact, and suppose that $X= X_1 + X_2$, with ${\mathcal N}$ being a dense subset of $X_2$.  Suppose also that $\phi: X \setminus \{\mathbf{1}\} \ra X$ is a   homeomorphism (with $\phi (\mathbf{n+1}) = \mathbf{n}$ for $\mathbf{n} \in \mathcal{N}$). Let $T: C (X) \ra C (X)$ be a codimension $1$ linear isometry such that $(Tf) (x) = f(\phi(x))$ for every $f \in C (X) $ and $x \in X \setminus \{\mathbf{1}\}$. Then $T$ is not an isometric shift.
\end{prop}

\begin{proof}
Since we are assuming that $T$ is a linear isometry, then there exists $\Delta $ in the dual space $C (X) '$, $\vc \Delta \vd \le 1$, such that $(Tf) (\mathbf{1}) = \Delta (f)$ for every $f \in C (X)$.
Take $\alpha_0 := \Delta( \xi_{X_1}) $ and $\alpha_1:= \Delta (\xi_{X_2})$, where $\xi_{X_i}$ denotes the characteristic function on $X_i$. If $\alpha_0 =0$, then define $f_0 := \xi_{X_1}$. Otherwise, there exists $\gamma \in \mathbb{K}$ such that $\gamma \alpha_0 + \alpha_1 =1$, and define $f_0 := \gamma \xi_{X_1} + \xi_{X_2}$. Then it is easy to check that in either case $Tf_0 =f_0$, so  $T$ is not a shift.
\end{proof}

As a consequence, we  see that the above proposition gives us a clue on  how to proceed in order to find isometric shifts of type $\mathrm{I}$ with $a=1$. This will be done next.

\subsection{A remark about notation}\label{pal}
In this section, we
deal with a compact space $Y$, for which we assume that there exists a   homeomorphism $\chi: Y \ra Y$ such that  the following two properties hold: 1)  It is  transitive with respect to a point $1_{Y} \in Y$, that is, the set $\tl \chi^n (1_{Y}) : n \in \mathbb{N} \tr$ is dense in $Y$; 2) there exists a  point $0_{Y} \in Y \setminus \{1_Y\}$ which is periodic for $\chi$, of prime period $N \in \mathbb{N}$. 

\medskip 
Spaces satisfying both requirements can be found for instance in  Example~\ref{delon}. On the one hand, we have   all separable infinite powers of any compact space with more than one point (see Corollary~\ref{purusa}). On the other hand, related to this,  chaotic homeomorphisms (in the Devaney sense, meaning in particular that they are transitive and have a dense set of periodic points; see \cite{BBCDS}) provide a rich source of examples. It is well known that, for $n \ge 2$,  every $n$-dimensional compact manifold admits a  homeomorphism of this kind (as proved in 
\cite{AD}), and indeed the set of all chaotic measure-preserving homeomorphisms on it is dense in the space of all measure-preserving homeomorphisms (see also \cite{DF, AS1}, and \cite{AS}). Finally, as announced in \cite[Remark 1]{AO}, these results give us a very different example satisfying the above requirements, as is the Sierpi\'nski curve.

\medskip

We denote $0_Y^0 := 0_Y, \ldots, 0_Y^{N-1} := \chi^{N-1} (0_Y)$. Next we take a compactification of $\mathcal{N}$ by $N$ points $\pmb{\infty}_0, \ldots, \pmb{\infty}_{N-1}$. Namely, for $k=0, \ldots, N-1$, define $\mathcal{N}_k := \{\mathbf{n} \in \mathcal{N} : n = k \mod N \}$, and consider the one-point compactification $\mathcal{N}_k \cup \{\pmb{\infty}_k\}$ of $\mathcal{N}_k$.

 Let $X$ be the union of $Y$ and ${\mathcal N} \cup \{\pmb{\infty}_0, \ldots, \pmb{\infty}_{N-1}\}$,  where each point $\pmb{\infty}_k$   is identified with the point $0_{Y}^k$ of $Y$. From now on we will just use the notation $\pmb{\infty}_k$ to denote this point both as a member of $Y$ and of ${\mathcal N} \cup \{\pmb{\infty}_0, \ldots, \pmb{\infty}_{N-1}\}$.

Define the map $\phi: Y \ra Y$ as $\phi := \chi^{-1}$.
Next take any two numbers $\delta_1, \delta_2 \in \mathbb{K}\setminus \{0\}$ satisfying $\va \delta_1 \vb + \va \delta_2 \vb \le 1$ and $(\delta_1 + \delta_2)^N \neq 1$, and define
 $\Delta : C (X)  \ra \mathbb{K}$  as 
\[ \Delta (f)  :=  \delta_1 f(1_{Y}) + \delta_2 f(\chi^N (1_{Y}))\]   for every $f$. It is immediate to see that $\Delta$ is linear and continuous, and that $ \vc \Delta  \vd  \le 1$.

Finally we define $T: C (X) \ra C (X)$ as $T:= T[1, \phi, \Delta]$. It is clear that $T$ is not an isometric shift of type $\mathrm{II}$ since we are assuming that both  $\delta_i$ are different from $0$.

\subsection{The results}

\begin{thm}\label{graciass0}
Assume that we follow the notation given in Subsection~\ref{pal}. Then $T$ is an isometric shift  of type $\mathrm{I}_0$.
\end{thm}

\begin{proof} 
Suppose that $f \in \bigcap_{n=1}^{\infty} T^n (C(X))$ satisfies $\vc f \vd \le 1$.
It is clear that for every $\mathbf{n} \in \mathcal{N}$, 
\begin{eqnarray}\label{rodillano}
f (\mathbf{n}) &=&(T^{-n+1} f ) (\mathbf{1})  \nonumber \\&=&  \delta_1 (T^{-n} f ) (1_{Y}) + \delta_2 (T^{-n} f ) (\chi^N (1_{Y})) \nonumber
\\ &=& \delta_1 f ( \chi^{n} (1_{Y})) + \delta_2 f ( \chi^{n+N} (1_{Y})).
\end{eqnarray}

Let us now see that $f (\pmb{\infty}_k) = 0$ for $k =0, \ldots, N-1$. For each $k$, define $\mathbb{N}_k := \{n \in \mathbb{N} : n = k \mod N \}$, and $\mathscr{M}_k := \{\chi^n (1_Y) : n \in \mathbb{N}_k\}$.
 We know that there exists $\lambda (0) \in \{0, \ldots, N-1\}$ such  that $\pmb{\infty}_0 $ belongs to the closure of $\mathscr{M}_{\lambda(0)}$.
Taking  into account that $\chi$ is  transitive and that $\chi^N (\pmb{\infty}_0) = \pmb{\infty}_0$, we deduce that 
there exists an increasing sequence $(n_i)$ in $\mathbb{N}_{\lambda(0)}$ such that both  $\pl f(\chi^{n_i} (1_{Y})) \pr $ and 
$\pl \pl f \circ \chi^N \pr  (\chi^{n_i } (1_{Y})) \pr $ converge to $f (\pmb{\infty}_0)$.  As a consequence,  we have  by Equality~\ref{rodillano} that
\begin{eqnarray*}
f (\pmb{\infty}_{\lambda(0)}) &=& \lim_{i \ra \infty} f (\mathbf{n_i}) \\
&=& \lim_{i \ra \infty} \delta_1 f ( \chi^{n_i} (1_{Y})) + \delta_2 f ( \chi^{n_i+N} (1_{Y})) \\
&=& (\delta_1 + \delta_2) f (\pmb{\infty}_0) .
\end{eqnarray*}
In the same way, we see that, if $N>1$, then for any $k \in \{1, \ldots, N-1\}$, there exists $\lambda (k)$ with 
$f (\pmb{\infty}_{\lambda(k)}) = (\delta_1 + \delta_2) f (\pmb{\infty}_k)$, and that $\lambda (k)$ can be taken equal to $\lambda (0) + k \mod N$. Then $f (\pmb{\infty}_{k}) = (\delta_1 + \delta_2)^N f (\pmb{\infty}_k)$ for every $k$,
which implies that $f (\pmb{\infty}_k) =0$.

Our next step consists of proving that $\delta_1 f(x) + \delta_2 f(\chi^N (x)) =0$ for every $x \in Y$. Notice first that $\lim_{n \ra \infty} f (\mathbf{n}) =0$, so by Equality~\ref{rodillano}, $\lim_{n \ra \infty} \delta_1 f ( \chi^{n} (1_{Y})) + \delta_2 f ( \chi^{n+N} (1_{Y})) =0$. Now, as above, if $x \in Y$, then there exists an increasing sequence $(m_i)$ in $\mathbb{N}$ such that $\pl f(\chi^{m_i} (1_{Y})) \pr $ and $\pl f(\chi^{m_i +N} (1_{Y})) \pr $ converge to $f (x)$ and $f (\chi^N (x))$, respectively. Consequently,
\begin{eqnarray*}
\delta_1 f(x) + \delta_2 f(\chi^N (x)) &=& \lim_{i \ra \infty}  \delta_1 f(\chi^{m_i} (1_{Y})) + \delta_2 f(\chi^{m_i +N} (1_{Y})) \\ 
&=& 0,
\end{eqnarray*}
as we wanted to see.

In particular, the above implies by Equality~\ref{rodillano} that $f (\mathbf{n}) = 0$ for every $\mathbf{n} \in \mathcal{N}$. It also implies that, for every $x \in Y$, 
$$f (x) = - \frac{\delta_2}{\delta_1}  f (\chi^{N} (x)),$$  so 
$$f (x) = \pl -1 \pr^n \pl \frac{\delta_2}{\delta_1} \pr^{n} f (\chi^{nN}(x))$$ for every $n \in \mathbb{N}$.
Obviously, to keep $f$ bounded, we need either that $f \equiv 0$ on $Y$  or that $\va \delta_1 \vb = \va \delta_2 \vb$.
But remark that in the latter case $\va f (x) \vb = \va f (\chi^{nN} (x)) \vb$ for every $x \in Y$ and 
every $n \in \mathbb{N}$. Now fix $k \in \{0, \ldots , N-1\}$, and notice that $\pmb{\infty}_k$ belongs to the closure of $\mathscr{M}_{\lambda(k)}$.  We deduce that $\va f (\pmb{\infty}_k) \vb = \va f (\chi^{\lambda(k)} (1_{Y})) \vb $, and that $f(\chi^{\lambda(k) +nN} (1_{Y})) = 0$ for every $n \in \mathbb{N} \cup \{0\}$. A similar process allows us to show that $f(\chi^{i+nN} (1_{Y})) = 0$ for every $i \in \{0, \ldots, N-1\}$ and $n \in \mathbb{N} \cup \{0\}$. By denseness, $f \equiv 0$ on $Y$, as we wanted to prove.
\end{proof}

The proof of the following corollary is now easy.

\begin{cor}\label{labaderu}
Let $n \in \mathbb{N}$. Let $Y_0$ be a connected and compact space with more than one point, and suppose that $\phi: Y_0 \ra Y_0$ is a  homeomorphism having a periodic point,  and such that $\phi^n$ is transitive. Then  there exist a compact space $X$ and an isometric shift  of type $\mathrm{I}_0$ on $ C(X)$ with $a \equiv 1$, such that $X \setminus \mathcal{N}$ consists exactly of $n$ connected components with more than one point, each homeomorphic to $Y_0$.
\end{cor}

\section{The case of the Cantor set}\label{niales}

In this section, we use results given in Sections~\ref{majorca} and \ref{luzia} to prove Theorem~\ref{toallas}.

\begin{proof}[Proof of Theorem~\ref{toallas}]
Without loss of generality we assume that $(x_n)$ is strictly monotone. Suppose first that $L \notin \mathbf{K}$. Then the result is an immediate consequence of Corollary~\ref{lanudo}, whose notation we follow.  If we take $E= \mathbb{R}$ and $A= \mathbf{K}$,  then $X_C$ will be 
a compact, metrizable and zerodimensional set (that is, having a basis consisting of closed and open sets), without isolated points. This implies that $X_C$ is itself homeomorphic to the Cantor set, and  we are done.

Suppose next that $L \in \mathbf{K}$. The realization of the Cantor set we will work with will be the ring $\mathbb{Z}_p$ of $p$-adic integers, for any $p \ge 2$,   which  is compact and zerodimensional when 
endowed with the metric given by its usual absolute value.
 We take  a  homeomorphism  $\mathbf{i} : \mathbb{Z}_p \ra \mathbf{K}$. By Corollary~\ref{purusa}, there are a homeomorphism 
$\phi: \mathbf{K} \ra \mathbf{K}$ and points $w, M$ in $\mathbf{K}$ such that $(\mathbf{K}, \phi, w)$ is totally transitive and $\phi (M) =M$. Let $L_0 := \mathbf{i}^{-1} (L), M_0 := \mathbf{i}^{-1} (M) \in \mathbb{Z}_p$, and define
$\mathbf{s}: \mathbb{Z}_p \ra \mathbb{Z}_p $ as $\mathbf{s} (t) : = t-L_0 +M_0$ for every $t \in \mathbb{Z}_p$,
which is obviously a homeomorphism. Moreover $\mathbf{j} := \mathbf{i} \circ \mathbf{s} \circ \mathbf{i}^{-1}$ is a homeomorphism from $\mathbf{K}$ onto itself such that $\mathbf{j} (L) = M$. As in  Lemma~\ref{aqua}, it is easy to see
that $(\mathbf{K},  \mathbf{j}^{-1} \circ \phi \circ \mathbf{j}, \mathbf{j}^{-1} (w))$ is totally transitive, and that
$\mathbf{j}^{-1} \circ \phi \circ \mathbf{j} (L) = L$.

Our final step consists of applying Theorem~\ref{graciass0} (where we identify $\mathbf{K}$ with $Y$ and $L$ with $0_Y = \pmb{\infty}$), and the proof is finished.
\end{proof}

\begin{rem}\label{exordo}
In Theorem~\ref{toallas}, when $L \notin \mathbf{K}$, there is no isometric shift with $a \equiv 1$. This is an immediate consequence of Proposition~\ref{nde}.
\end{rem}

\section{Acknowledgements}
The author would like  to thank Professor Francisco Santos for useful conversations.

\end{document}